\documentclass[a4paper,twoside,12pt,fleqn]{scrartcl}

\usepackage[italian,ngerman,english]{babel}       
\usepackage[utf8]{inputenc}	
\usepackage[T1]{fontenc}

\usepackage[english,cleanlook]{isodate}

\usepackage{lmodern} 
\usepackage[babel=true]{microtype} 

\usepackage{hyperref} 
\hypersetup{
hidelinks 
}

\usepackage{graphicx} 
\usepackage{tikz}

\usepackage{enumitem}

\RedeclareSectionCommand[
	afterskip=-1em, 
]{subsubsection}

\usepackage{mathtools} 
\usepackage{amsthm} 
\usepackage{amssymb} 

\usepackage{bbm} 
\usepackage{mathrsfs} 

\newcommand{\R}{\mathbb{R}}
\newcommand{\C}{\mathbb{C}}
\newcommand{\N}{\mathbb{N}}

\newcommand{\e}{\mathrm{e}}  

\DeclareMathOperator{\divop}{div}

\DeclareMathOperator{\dist}{dist}

\DeclareMathOperator{\id}{id}

\DeclareMathOperator{\tr}{tr}

\DeclareMathOperator{\Aut}{Aut}

\newcommand{\loc}{\text{loc}}

\def\Xint#1{\mathchoice
        {\XXint\displaystyle\textstyle{#1}}%
        {\XXint\textstyle\scriptstyle{#1}}%
        {\XXint\scriptstyle\scriptscriptstyle{#1}}%
        {\XXint\scriptscriptstyle\scriptscriptstyle{#1}}%
        \!\int}
\def\XXint#1#2#3{{\setbox0=\hbox{$#1{#2#3}{\int}$}
\vcenter{\hbox{$#2#3$}}\kern-.5\wd0}}
\def\dashint{\Xint-}

\numberwithin{equation}{section}

\theoremstyle{plain}
\newtheorem{definition}{Definition}[section]
\newtheorem*{definition*}{Definition} 

\newtheorem{theorem}[definition]{Theorem}
\newtheorem*{theorem*}{Theorem} 
\newtheorem{lemma}[definition]{Lemma}
\newtheorem{corollary}[definition]{Corollary}
\newtheorem{proposition}[definition]{Proposition}
\newtheorem{definition-proposition}[definition]{Definition--Proposition}

\newtheorem*{question*}{Question}
\newtheoremstyle{boldremark} 
  {3pt}
  {3pt}
  {}
  {}
  {\bfseries}
  {}
  {.5em}
  {}
\theoremstyle{boldremark}

\newtheorem{remark}[definition]{Remark}

\newenvironment{bfproof}[1][\proofname]{
	\begin{proof}[\normalfont\bfseries #1]}{\end{proof}
	} 

\usepackage[left=2cm, right=2cm, bottom=3cm, top=3cm]{geometry}

\title{
The Parametric Approach\\
to the Willmore Flow}

\date{}

\author{
Francesco Palmurella
\and
Tristan Rivi\`ere}

\usepackage[headsepline]{scrlayer-scrpage}
\clearscrheadfoot 
\automark[section]{section} 
\ohead{\textup\headmark}
\usepackage{duerer}

\mathtoolsset{showonlyrefs}

\begin{document}

\maketitle

\begin{abstract}
We introduce a parametric framework
for the study of Willmore gradient flows
which enables to consider a general class of
weak, energy-level solutions and opens
the possibility to study energy quantization
and finite-time singularities.
We restrict in this first work to a small-energy
regime and prove that,
for small-energy weak immersions,
the Cauchy problem
in this class admits a unique solution.
\end{abstract}
\begin{center}
\medskip
\end{center}

{\small

\setcounter{tocdepth}{1}
\tableofcontents
}

\newpage
\section{Introduction}

The present paper moves the first steps
towards a parametric theory for the  Willmore flow that,
we believe, will lead to an effective study of singularities,
bubbling analysis and energy quantization.

We start by
recalling what a Willmore surface is.

\subsubsection{Willmore Surfaces}
The \emph{Willmore energy}
of a surface $\mathcal{S}\subset\R^n$
was
initially considered in the work of 
\textsc{Poisson}
\cite{PoissonMemoireElastique}
and \textsc{Germain} \cite{GermainReCherchesElastiques}
on elastic plates.
It was reconsidered in a purely geometric perspective
in the works of \textsc{Thomsen} \cite{ThomsenWillmore1923} and \textsc{Blaschke} \cite{BlaschkeVorlDiffIII},
in attempt to merge the study of minimal surfaces
with conformal invariance.
It was reintroduced in recent times
in the works of \textsc{Willmore} \cite{MR0202066},
and \textsc{Bryant} \cite{MR772125}
in pure mathematics
and by \textsc{Canham} \cite{CANHAM1970}
and \textsc{Helfrich} \cite{HELFRICH1973}
in theoretical biology while modeling the shapes of blood cells 
and lipid bilayers
(and, in fact, it is sometimes referred to as Canham-Helfrich energy).
It comes in the following three variants:
\begin{align}
\label{eq:Willmore-energy}
\mathcal{W}_0(\mathcal{S})
&=\frac{1}{2}\int_{\mathcal{S}} |A^\circ|^2\,d\sigma,&
\mathcal{W}_1(\mathcal{S})
&=\int_{\mathcal{S}} |H|^2\,d\sigma, &
\mathcal{W}_2(\mathcal{S})
&=\frac{1}{4}\int_{\mathcal{S}} |A|^2\,d\sigma,
\end{align}
were $H$ the mean curvature,
 $A$ is the 2nd fundamental form,
$A^\circ$ it tracefree part
and $d\sigma$ the area element.
If $K$ denotes the Gauss curvature of $\mathcal{S}$,
there holds
\begin{align}
\label{eq:pointwise-relations-A}
\frac{1}{2}|A^\circ|^2 
=|H|^2 - K
=\frac{1}{4}|A|^2 - \frac{1}{2}K,
\end{align}
hence, by Gauss-Bonnet theorem, if the topology
of $\mathcal{S}$ is fixed, at least in a smooth setting
such energies are all variationally equivalent
(in particular, they have the same Euler-Lagrange operator).
Depending on the context,
it may however be more favourable to work with one than another.

The \emph{Willmore operator} is
the associated Euler-Lagrange operator:
\begin{align}
\label{eq:Willmore-equation}
\delta\mathcal{W}=\Delta^\perp H + Q(A^\circ)H,
\end{align}
where $\Delta^\perp$ is the Laplace  operator on the normal
bundle and
\begin{align}
 Q(A^\circ)H
=\big\langle A^\circ, \langle H,A^\circ\rangle\big\rangle
= g^{\mu\sigma}g^{\nu\tau}\big\langle
A^\circ_{\mu\nu}, \langle A^\circ_{\sigma\tau},H\rangle
\big\rangle.
\end{align}
Similarly as for the mean curvature,
$\delta{\mathcal{W}}$ is a normal-valued vector field
along $\mathcal{S}$.
When $n=3$,
the expression simplifies somewhat:
\begin{align}
\label{eq:Willmore-equation-codim-1}
\delta\mathcal{W}
=\Delta^\perp H + |A^\circ|^2H
=(\Delta H_{\text{sc}} + 2(H_{\text{sc}}^2-K)H_{\text{sc}})N,
\end{align}
where
$N$ is the Gauss map of $\mathcal{S}$ and
$H_{\text{sc}}=\langle H,N\rangle$ is
the scalar mean curvature.

\emph{Willmore surfaces} are
those surfaces with vanishing Willmore operator.
Any of the Willmore energies \eqref{eq:Willmore-energy}
is invariant under conformal transformations
of $\R^n$ and,
in fact, the Lagrangian density
$|A^\circ|^2d\sigma$ is a pointwise conformal invariant,
see \textsc{Chen} \cite{MR0370436}.
As a consequence, the Willmore operator and the
 notion of Willmore surface are also conformal invariants.
\medskip

\subsubsection{Willmore Flows}
A \emph{Willmore $L^2$--gradient flow}
in $\R^n$
(\emph{Willmore flow} for short) of a
closed, abstract surface $\Sigma$
is a
1-parameter family of immersions
$\Phi(t,\cdot):\Sigma\to \R^n$, $t\in I\subseteq\R$ evolving
according to the law
\begin{align}
\label{eq:Willmore-flow}
\frac{\partial}{\partial t}\Phi
=-\delta\mathcal{W}
+U
\quad\text{in } I\times \Sigma,
\end{align}
where, for each $t$, $\delta\mathcal{W}$ is the
Willmore operator of $\mathcal{S}_t=\Phi(t,\Sigma)$
and $U=U^\mu\partial_\mu\Phi$ is a
tangent tangent vector field,
possibly time-dependent.
\medskip

One good reason to consider Willmore flows is
that they satisfy the energy identity, namely
if $I=(0,T)$, then
\begin{align}
\label{eq:energy-identity}
\mathcal{W}_0(\Phi(t,\cdot))
-\mathcal{W}_0(\Phi(0))
=-\int_0^t\int_{\Sigma}
|\delta\mathcal{W}|^2\,d\sigma_g\,d\tau,
\quad\text{for } 0\leq t < T.
\end{align}
It may be rephrased by saying that, 
among all families of immersions whose velocity vector 
has normal part with
$L^2$-norm equal to $\|\delta\mathcal{W}\|_{{L^2}(S^2)}$,
Willmore flows are those with
most rapidly decreasing Willmore energy
(in any of the forms given in \eqref{eq:Willmore-energy}).
Thus, at least in principle,
they have the potential to converge
efficiently to Willmore immersions as $t\to+\infty$.

This is a feature common to gradient
flows that makes them  particularly worth studying.
The first to consider  $L^2$--gradient
flows in a geometric context were 
\textsc{Eells} and \textsc{Sampson} \cite{MR164306}
in the context of harmonic maps.
Since then, the study of parabolic geometric flows has widened
to the extent that some of them constitute research areas
on their own right, the
mean curvature flow and Hamilton's Ricci flow being
two of the best-known examples.
\medskip

It should be noted right away that what
is typically called a Willmore flow
is a family solving \eqref{eq:Willmore-flow}
with $U=0$,
which we will call here a 
\emph{normal Willmore flow}.
Since $\Sigma$ is closed, and $\delta\mathcal{W}$
is a tensor, it is classical fact that there is a bijective correspondence between tangential components
and family of reparametrizations of $\Sigma$,
see e.g. \textsc{Mantegazza} \cite[Proposition 1.3.4]{MR2815949} for the case,
entirely analogous in this regard, of the mean curvature flow.
Consequently, if, say, $I$ is a connected interval containing $0$,
for every family solving \eqref{eq:Willmore-flow}
there is a unique family of diffeomorphisms $\varphi:I\times\Sigma\to\Sigma$
with $\varphi(0,\cdot)=\id_{\Sigma}$ so that 
the reparametrized family $\Phi(t,\varphi(t,\cdot))$, $t\in I$ 
is a normal Willmore flow,
and on the other hand, every reparametrization
of a normal Willmore flow will be
a Willmore flow \eqref{eq:Willmore-flow}
for some $U$.

Thus, in this sense, similarly
as for immersions of surfaces, flows can be regarded
as equivalence classes of solutions to
\eqref{eq:Willmore-flow},
two of them being equivalent if one
can be reparametrized into another.
As for surfaces, depending on the 
situation one may choose one parametrization
over another, and in this case
this may be done through the choice of the tangential
component.
This will be a crucial fact in the present work.
\medskip

\medskip

The study of Willmore flows was introduced 
by \textsc{Kuwert} and \textsc{Sch\"atzle}
\cite{MR1882663,MR1900754}
and \textsc{Simonett} \cite{MR1827100}
and is since then subject of a growing number of
works.
Particularly useful for us will be the
one by \textsc{Kuwert} and \textsc{Scheuer}
\cite{10.1093/imrn/rnaa015}
providing asymptotic estimates on
the barycenter along the flow.

Our attention here focuses on the following
foundational result.
Consider the Cauchy problem for the normal Willmore flow:
\begin{align}
\label{eq:Cauchy-normal-Willmore}
\left\{
\begin{aligned}
\frac{\partial}{\partial t}\Phi & = -\delta\mathcal{W},
&&\text{in }(0,T)\times \Sigma,\\
\Phi(0,\cdot) &= \Phi_0 
&&\text{on }\Sigma.\\
\end{aligned}
\right.
\end{align}
\begin{theorem*}[\cite{MR1882663,MR1900754}]
There exists $\varepsilon_0(n)>0$
so that, for a smooth immersion ${\Phi_0:\Sigma\to\R^n}$
satisfying $\mathcal{W}_0(\Phi_0)=\mathcal{W}_0(\Phi_0(\Sigma))<\varepsilon_0$,
then \eqref{eq:Cauchy-normal-Willmore}
has a unique solution in the smooth category,
which furthermore exists for all times and converges
to a round sphere.
\end{theorem*}
It should be said immediately that
if $\mathcal{W}_0(\Phi)$ is sufficiently small, 
 $\Sigma$ must be a sphere.
Indeed, as already noticed
in \cite{MR0202066},
it is always $\mathcal{W}_1(\Phi)\geq 4\pi$
and from \eqref{eq:pointwise-relations-A},
using Gauss-Bonnet one sees that
\begin{align}
\label{eq:consequence-topology}
\chi(\Sigma)
=\frac{1}{2\pi}\int_{\Sigma} K\,d\sigma
=\frac{1}{2\pi}\big(\mathcal{W}_1(\Phi)-\mathcal{W}_0(\Phi)\big)
\geq\frac{1}{2\pi}
\big(4\pi-\varepsilon_0\big)
>1,
\end{align}
where $\chi(\Sigma)$ is the Euler-Poincar\'e 
characteristic of $\Sigma$.
\medskip

Such theorem was concerned with
\emph{smooth} solutions.
As for other geometric flows however,
an effective study of singularities
and bubbling analysis
requires eventually to work at the \textit{energy level},
namely, to consider appropriate notions
of weak solution.

We have in mind as a particular example
the classical work on
the harmonic map flow
done by \textsc{Struwe} \cite{MR826871,MR2431434}
and complemented by the works of others,
the 2nd author \cite{RiviereThesisIV}, \textsc{Freire} \cite{MR1384838,MR1324632},
\textsc{Chang}, \textsc{Ding}
and \textsc{Ye} \cite{MR1180392},
\textsc{Topping} \cite{MR1883901}
\textsc{Bertsch},
\textsc{Dal Passo}
and \textsc{van der Hout} \cite{MR1870959} just to mention a few.
\medskip

We believe that the framework introduced
by the 2nd author in a series
of works \cite{MR2430975,MR3276154,MR3524220},
which led for instance
to an effective energy quantization analysis
of Willmore surfaces by \textsc{Bernard} and the 2nd author
\cite{MR3194812} to be,
when suitably adapted, the appropriate one.
We want to give in the present paper an idea
of why this should  be true by introducing,
under
particularly favourable hypotheses,
an energy-level class of weak Willmore flow
and prove a uniqueness statement for the corresponding
Cauchy problem in this class
for a broad set of weak initial data,
which we believe to be sufficiently close to the
largest possible one
(among unbranched surfaces).

\medskip
Let us mention that
\textsc{Lamm} and \textsc{Koch}
in 
\cite{MR2916362} 
obtained (among other results of geometric interest)
an existence and uniqueness result
for the Willmore flow for entire graphs
in a weak framework  with
Lipschitz initial datum.
Such datum needs to be small in the
Lipschitz norm.

\subsubsection{Well-Balanced Conformal Willmore Flows}
We shall work, in the present paper,
always in a low energy regime,
namely we shall arrange things so that the Willmore energy
of the surfaces in consideration
$\mathcal{W}_0(\mathcal{S})$ is as small as needed;
furthermore, we shall also work in codimension one, namely $n=3$.
The first major consequence of this is that,
as already said above, with \eqref{eq:consequence-topology}
we may directly assume that the underlying
topology is that of the standard sphere $S^2$.
The second one is that we can take advantange of 
results from the work of 
\textsc{De Lellis} and \textsc{M\"uller}
\cite{MR2169583,MR2232206}.

So, from now, it is $\Sigma=S^2$,
and the underling reference metric and
complex structure are the standard ones.
\medskip

Central in the theory developed in 
\cite{MR2430975,MR3276154,MR3524220}
and in the present one
is the idea of working with \emph{conformal
immersions}.
The first advantage of
doing so is that the Willmore operator
becomes uniformly elliptic,
with ellipticity constants depending on the conformal factor,
and this permits eventually the regularity bootstrap.
The second one is that, 
exploiting conservation laws issuing from 
the conformal
invariance  (as explained in \textsc{Bernard} \cite{MR3518329}),
the Willmore operator
of a conformal immersion, which is a
4th order quasilinear elliptic system,
can be recast as a 2nd order semilinear system involving Jacobian-type
nonlinearities, which allows regularity bootstrap
by means of integrability by compensation,
similarly as in the work of \textsc{H\'elein} \cite{MR1913803} on 
weakly harmonic maps in two dimensions.

The idea is then to
consider 
\textit{Willmore flows in conformal gauge},
where the equation becomes uniformly parabolic,
if the conformal factor is uniformly bounded away from zero,
 and then use a slice-wise in time (elliptic)
integrability by compensation arguments to 
bootstrap the regularity of the equation,
which -- as is often the case when working with parabolic PDEs 
in small energy regime -- will suffice to get the regularity
also in the time variable.
This approach was successfully used by
the 2nd author in \cite{RiviereThesisIV} 
for the case of the harmonic map flow.
\medskip

Indeed,
\eqref{eq:Willmore-flow} is invariant under reparametrizations
and thus it is degenerate parabolic, as is the case for others geometric flows
such as the mean curvature flow or the Ricci flow - 
and this can be a serious source of troubles.
The celebrated trick of \textsc{DeTurck}
\cite{MR697987}, originally 
devised for the Ricci flow 
but easily adapted to the present situation,
is one way of overcoming this problem.
For the sake of complenetess, we outlined it in
Appendix \ref{sec:DeTurck}.
Such method
has the advantage of working
regardless of the topology of $\Sigma$,
but, as an inspection of the proof reveals,
does not seem to be suitable
when working with low degrees of smoothness
and moreover
does not give explicitly a control on the
parametrization that is chosen by such gauge.
\medskip
 
We are instead going to consider the Cauchy problem
\begin{align}
\label{eq:Cauchy-conformal-Willmore}
\left\{
\begin{aligned}
\frac{\partial}{\partial t}\Phi & = -\delta\mathcal{W}+U,
&&\text{in }(0,T)\times S^2,\\
\Phi(0,\cdot) &= \Phi_0 
&&\text{on }S^2.\\
\end{aligned}
\right.
\end{align}
where the tangential vector field
$U$ is chosen so that $\Phi(t,\cdot)$
is conformal for every $t$.
\medskip

It is a fairly simple matter to find
an explicit characterization for $U$,
which --
 unsurprisingly, given the relationship between complex and conformal
structures on surfaces
-- is best expressed in complex notation.
To this aim, we recall that
if $\Phi:S^2\to\R^3$ is a conformal immersion
with metric $g=\e^{2\lambda}g_{S^2}$,
the second fundamental form
may be written in complex notation as
\begin{align}
A
=h_0 + \overline{h}_0 +
H\otimes g,
\end{align}
where $H$ is the mean curvature and
\begin{align}
h_0 = A_{zz}\,dz\otimes dz
=\frac{1}{4}(A_{11} - A_{22} - 2iA_{12})dz\otimes dz.
\end{align}
Similarly, the tracefree second fundamental form is written as
\begin{align}
A^\circ = A - Hg = h_0 + \overline{h}_0.
\end{align}
With this formalism, we have the following.
\begin{lemma}
\label{lemma:conformal-flow-cond}
The tangential component of a conformal Willmore flow
satisfies
\begin{align}
\label{eq:conformal-flow-cond}
\bar{\partial}U^{(1,0)} = 
-\langle
\delta\mathcal{W}, \overline{h}_0{}
\rangle^{\sharp_g},
\end{align}
that is
\begin{align}
\label{eq:conformal-flow-cond-loc}
\partial_{\bar{z}} (U^1+iU^2)\partial_z\otimes d\bar{z}
=2\e^{-2\lambda}\big\langle
-\delta\mathcal{W}, A_{\bar{z}\bar{z}}
\big\rangle\partial_z\otimes d\bar{z}.
\end{align}
\end{lemma}
Basic facts for the
$\overline{\partial}$-operator on vector fields
is recalled in Appendix \ref{sec:delbar}.
It is certainly good news that,
at least on the sphere,
it defines an uniformly elliptic,
zero-cokernel operator.
\medskip

However,
\eqref{eq:conformal-flow-cond} do not suffices \textit{per se}
to guarantee that
\eqref{eq:Cauchy-conformal-Willmore} is uniformly parabolic,
not even for short time,
since the control on the conformal factor 
in time still needs to be addressed.
Geometrically this is evident:
any conformal Willmore flow
can be composed with any 1-parameter family
of conformal self-maps of $S^2$ and remain conformal,
and $\Aut(S^2)$, the set conformal
(i.e. biholomorphic) self-map of
$S^2$ is not compact.
One can clearly see this also by noting that 
\eqref{eq:conformal-flow-cond}
has a nontrivial kernel 
given by the conformal Killing
(holomorphic, with complex identification) vector fields
on $S^2$
which is 6 dimensional,
as $\Aut(S^2)$.

Moreover, 
\eqref{eq:conformal-flow-cond}
is an equation satisfied for every fixed $t$,
and does not give any information on the
regularity in time of $U$.
Geometrically, this this means that
the 1-parameter family
of maps in $\Aut(S^2)$ which we may compose a conformal 
Willmore flow may be taken nonsmooth with respect to $t$.
\medskip

Precisely because $\Aut(S^2)$ is 6-dimensional
however,
a final choice of a 6-dimensional constraint
will be enough to tame the action of such gauge group.
This will be defined by the following.
\begin{definition}
\label{def:well-balanced}
An immersion $\Phi:S^2\to\R^3$
is called 
\emph{well-balanced}
if there holds
\begin{align}
\label{eq:well-balanced}
\int_{S^2} I d\sigma_{g}=0
\quad\text{and}\quad
\int_{S^2} \Phi\times I\,d\sigma=0,
\end{align}
where $I$ denotes  standard embedding of $S^2$, $d\sigma$ its
area element
and $d\sigma_g$ the area element
for the induced metric $g=\Phi^\star g_{\R^3}$.
\end{definition}

\begin{remark}
\label{rmk:balance-transl}
Note that being well-balanced is a translation-invariant condition,
namely if $\Phi$ is well-balanced, so is $\Phi + k$ for every
$k\in\R^3$.
\end{remark}

Conditions
\eqref{eq:conformal-flow-cond}
and \eqref{eq:well-balanced}
together with a good choice for the parametrization
of the initial datum , which we shall now discuss,
will be sufficient to control the
behaviour of  the tangential component $U$ in
on \eqref{eq:Cauchy-conformal-Willmore}.
Moreover,
they are meaningful also
for the notion of weak
conformal Willmore flow that we are going to define.

\subsubsection{
Chosing an ad-hoc Parametrization for initial Data with Small Energy.}
From a geometric perspective,
both the Cauchy problems \eqref{eq:Cauchy-normal-Willmore}
and \eqref{eq:Cauchy-conformal-Willmore}-\eqref{eq:conformal-flow-cond}
possess an obvious ``gauge invariance''
for the initial datum, namely
if $\Phi_0(S^2)=\mathcal{S}$ is the immersed sphere
representing the initial datum, and
$\varphi$ is any diffeomorphism of $S^2$,
then $\Phi_0\circ\varphi$ is
again a parametrization for the same surface $\mathcal{S}$,
and there is no a priori preferred choice -- or possibility
to distinguish -- between
$\Phi_0$ and  $\Phi_0\circ\varphi$.
This a relevant issue for a parametric theory.

The conformal gauge choice
helps to reduce this invariance
($\Phi_0$ has to be conformal,
and so $\varphi$ must belong to $\Aut(S^2)$),
but does not break it entirely.
To this aim, we shall use
the following result contained in the
work of 
\textsc{De Lellis} and \textsc{M\"uller}
\cite{MR2169583,MR2232206}.

\begin{theorem*}[\cite{MR2169583,MR2232206}]
There exist $\varepsilon_0,\,C\,>0$ so that, if
$\mathcal{S}\subset\R^3$ is an immersed surface
with area 
$\mathcal{A}(\mathcal{S})=4\pi$
and Willmore energy $\mathcal{W}_0(\mathcal{S})\leq\varepsilon_0$,
there exists a conformal parametrization
$\Phi:S^2\to\mathcal{S}$ satisfying
\begin{align}
\label{eq:DeLellis-Muller-intro}
\|\Phi - I - c\|_{W^{2,2}(S^2)}+
\|\e^{\lambda}-1\|_{L^\infty(S^2)}
\leq C\sqrt{\mathcal{W}_0(\mathcal{S})},
\end{align}
where $I:S^2\to \R^3$ denotes the standard immersion
of $S^2$ and $c=\dashint_{S^2}\Phi\,d\sigma$.
\end{theorem*}
In this theorem the fact that
the area of the surface is $4\pi$
can be seen to a normalization
achievable by scaling.
We shall need another one, achievable by translations;
to this aim recall that the \emph{barycenter} of an immersed surface
$\mathcal{S}\subset\R^3$ is defined as
\begin{align}
\mathcal{C}(\mathcal{S})
=\int_{\mathcal{S}} \id_{\R^3}\, d\mathcal{H}^2
=\int_{\Sigma}\Phi\, d\sigma_g,
\end{align}
where $\Phi:\Sigma\to\mathcal{S}$ is any parametrization of $\mathcal{S}$.

The set of initial data for 
the conformal Willmore flow will consist
\textit{geometrically} of the set of 
immersed surfaces $\mathcal{S}\subset \R^3$
with Willmore energy $\mathcal{W}_0(\mathcal{S})\leq\varepsilon$,
area $4\pi$ and vanishing barycenter $\mathcal{C}(\mathcal{S})=0$.
\textit{Parametrically} we shall choose
a parametrization provided by the above theorem,
which is in addition well-balanced as in
Definition \ref{def:well-balanced}.
More precisely:
\begin{definition}
\label{eq:D-class}
For $\varepsilon>0$, 
$\mathscr{D}^\varepsilon(S^2,\R^3)$
is the set of smooth conformal immersions
$\Phi:S^2\to\R^3$ so
that the surface $\mathcal{S}=\Phi(S^2)$ has
Willmore energy $\mathcal{W}_0(\mathcal{S})\leq \varepsilon$,
area $\mathcal{A}(\mathcal{S})=4\pi$,
barycenter $\mathcal{C}(\mathcal{S})=0$,
is well-balanced
and so that \eqref{eq:DeLellis-Muller-intro} holds
for $C>0$ given by that estimate.
\end{definition}

This is,
when restricted to the smooth category,
the suitable class
of initial data
that shall be considered in this work,
for sufficiently small $\varepsilon>0$.
It will be enlarged
to its weak $W^{2,2}$-closure
when considering the extension of the theory
to the weak framework,
which we discuss below.
\medskip

We want to stress that the only essential requirement
in Definition \ref{eq:D-class}
is the control (smallness) of the Willmore energy.
All the others can be seen as normalizations.
More precisely, the first result of this work,
which will be used to prove the main one, is 
the following
extension of the theorem above:
\begin{proposition}
\label{prop:well-balanced-fam}
There are $\varepsilon_0,\,\delta,\, C > 0$ 
with the following properties:
\begin{enumerate}[label=(\roman*)]
\item Any immersed surface $\mathcal{S}\subset\R^3$
with area $4\pi$ and Willmore energy 
$\mathcal{W}_0(\mathcal{S})\leq\varepsilon_0$
admits a conformal parametrization satisfying
\eqref{eq:DeLellis-Muller-intro}
which is also well-balanced.
\item 
For any well-balanced conformal
parametrization $\Psi: S^2\to\mathcal{S}$
with conformal factor $\e^\nu$ 
and some vector $c\in\R^3$ so that
\begin{align}
\label{eq:WB-1}
\|\Psi - I-c\|_{W^{2,2}(S^2)}
+\|\e^{\nu} -1\|_{L^\infty(S^2)}
\leq\delta,
\end{align}
there holds
\begin{align}
\label{eq:WB-final}
\|\Psi - I-c\|_{W^{2,2}(S^2)}
+\|\e^{\nu} -1\|_{L^\infty(S^2)}
\leq C\sqrt{\mathcal{W}_0(\mathcal{S})}.
\end{align}
Furthermore, the following \emph{local uniqueness}
property holds: if $\Psi'$ is another well-balanced 
conformal immersion satisfying \eqref{eq:WB-1},
and $\psi\in\Aut(S^2)$ is the conformal
diffeomorphism so that $\Psi'=\Psi\circ \psi$,
there is a neighborhood $\mathcal{O}\subset\Aut(S^2)$ of the identity
$e$
(depending only on $\delta$)
so that, if $\psi\in \mathcal{O}$, then $\psi=e$.
\end{enumerate}
\end{proposition}

\subsubsection{Conformal Weak Flows}
We now define an energy-level class
of maps where one can consider
 weak conformal Willmore flows.
We believe it to the a prototype for 
future works concerned with Willmore flows
at energy level.

Central to the definitions we shall give shortly
is that the validity of the
energy identity \eqref{eq:energy-identity}
 (in fact, a slightly weaker version will suffice).
This should be, broadly speaking,
a requirement to avoid the presence of 
pathological solutions that invalidate the uniqueness
of the solution to the Cauchy problem,
as the examples of \textsc{Topping} \cite{MR1883901}
and \textsc{Bertsch},
\textsc{Dal Passo}
and \textsc{van der Hout} \cite{MR1870959}
show in the case of the harmonic map flow.
\medskip

From \cite{MR2430975,MR3276154,MR3524220}
we recall the notion of
weak $W^{2,2}$-Lipschitz immersion.
If
$W^{1,\infty}_{\text{imm}}(S^2,\R^3)$ denotes
the set of Lipschitz immersions, namely those Lipschitz
maps $\Phi:S^2\to\R^3$ so that
there exists $C=C(\Phi)>0$ with
\begin{align}
\frac{1}{C}g_{S^2}
\leq g=\Phi^\star g_{\R^3}
\leq C g_{S^2}
\end{align}
almost everywhere in the sense of metrics,
we let
\begin{align}
\label{eq:E}
\mathscr{E}(S^2,\R^3)
=W^{1,\infty}_{\text{imm}}(S^2,\R^3)\cap W^{2,2}(S^2).
\end{align}
Every map in such set admits a conformal reparametrization
and moreover, it is possible to define its Willmore operator
in the sense of distributions.
Starting from the divergence form of the Willmore operator
introduced in \cite{MR2430975}:%
\footnote{
We denote, here and in  the sequel:
\begin{enumerate}[label=--=]
\item
$
\nabla^{*_g}(Z\otimes \omega)
= \frac{1}{\sqrt{ g}}\partial_\mu
\big(\sqrt{g}\,g^{\mu\nu}\omega_\nu Z\big),
$
\emph{minus} the formal $L^2$-adjoint of the covariant derivative
induced on the pull-back bundle $\Phi^\star(T\R^3)$
acting on sections of $\Phi^\star(T\R^3)\otimes T^*S^2$,
\item
$
\langle A,H\rangle^{\sharp_g}
=g^{\mu\xi}
\langle A_{\xi\nu},H\rangle\partial_{\mu}\otimes dx^\nu
\simeq g^{\mu\xi}
\langle A_{\xi\nu},H\rangle\partial_{\mu}\Phi\otimes dx^\nu$
the the 1st-index raising  of $\langle A,H\rangle$,
and similarly for $\langle A^\circ,H\rangle$.
\end{enumerate}
}
\begin{align}
\delta\mathcal{W}
= \nabla^{*_g}\big( 
\nabla H
-2(\nabla H)^{\top_\Phi}
-|H|^2d\Phi
\big)
=\nabla^{*_g}\big( \nabla H 
+\langle A^\circ,H\rangle^{\sharp_g}
+\langle A,H\rangle^{\sharp_g}\big),
\end{align}
the Willmore operator of $\Phi\in\mathscr{E}(S^2,\R^3)$
is defined as the distribution-valued two form
given by
\begin{align}
\label{eq:Weak-willmore}
\big(\delta\mathcal{W}d\sigma_g, \varphi
\big)_{\mathcal{D}'}
&=\int_{S^2}\Big(
\big\langle H,\Delta_g\varphi\big\rangle
-\big\langle\langle A^\circ, H\rangle^{\sharp_g}
+\langle A, H\rangle^{\sharp_g},\nabla\varphi\big\rangle_g
\Big)d\sigma_g,
\end{align}
for every $\varphi\in C^\infty(S^2,\R^3)$.
\medskip

We can now give the following
central definitions.
\begin{definition}[Weak Initial Data]
\label{def:weak-initial-data}
For $\varepsilon>0$, 
$\mathscr{W}^\varepsilon(S^2,\R^3)$
is the closure of
$\mathscr{D}^\varepsilon(S^2,\R^3)$
(from Definition \ref{eq:D-class})
with respect to the weak $W^{2,2}(S^2)$-topology.
\end{definition}

We will consider $\mathscr{W}^\varepsilon(S^2,\R^3)$
only for 
$\varepsilon>0$ sufficiently small.
As a consequence of the works
\cite{MR2430975,MR3276154,MR3524220},
$\mathscr{W}^\varepsilon(S^2,\R^3)$ is a subset of the space $\mathscr{E}(S^2,\R^3)$,
which would be the broadest possible choice
for the theory
(among nonbranched surfaces, at least).
We do not know at present whether $\mathscr{W}^\varepsilon(S^2,\R^3)$
coincides, or strictly contained in,
$\mathscr{E}(S^2,\R^3)$.

\begin{definition}[Well-Balanced Energy Class]
\label{def:energy-class}
For $\varepsilon,\,\delta,\,T>0$,
$\mathscr{W}_{[0,T]}^{\varepsilon,\delta}(S^2,\R^3)$
is set of locally integrable maps 
${\Phi:(0,T)\times S^2\to\R^3}$
so that
\begin{enumerate}[label=(\roman*)]
\item For almost every $t$,
$\Phi(t,\cdot)$ is in  $\mathscr{E}(S^2,\R^3)$
and conformal,
\item There holds
\begin{align}
\label{eq:def-class-DLM}
\|\Phi - I-c\|_{L^\infty((0,T), W^{2,2}(S^2))}
+\|\e^\lambda-1\|_{L^\infty((0,T)\times S^2)}\leq \delta,
\end{align}
where
$I$ denotes  standard embedding of $S^2$, $\e^{\lambda}=\e^{\lambda(t,\cdot)}$ is the
conformal factor of $\Phi(t,\cdot)$
and $c(t)=\dashint_{S^2}\Phi(t,\cdot)\,d\sigma$,

\item There holds
\begin{align}
\label{eq:weaker-energy-id}
\delta\mathcal{W}\in L^2((0,T)\times S^2)
\quad\text{and}\quad\mathcal{W}_0(\Phi(t,\cdot))\leq \varepsilon
\;\text{ for a.e. } t,
\end{align}
\item $\Phi$ is well-balanced for a.e. $t$.
\end{enumerate}
Finally we let, also for $T=+\infty$,
\begin{align}
\mathscr{W}_{[0,T)}^{\varepsilon,\delta}(S^2,\R^3)
=\bigcap_{\tau\in(0,T)}\mathscr{W}_{[0,\tau]}^{\varepsilon,\delta}(S^2,\R^3).
\end{align}
\end{definition}

Assumption \eqref{eq:def-class-DLM}
is quite natural if we look at Proposition 
\ref{prop:well-balanced-fam}.
In the energy class, a weak Willmore flow is defined 
as follows.

\begin{definition}
[Weak Willmore Flow]
\label{def:weak-Willmore-flow}
$\Phi\in\mathscr{W}_{[0,T)}^{\varepsilon,\delta}(S^2,\R^3)$ 
is a weak solution of the Willmore flow with tangential component
$U=U^\mu\partial_\mu\Phi$:
\begin{align}
\frac{\partial}{\partial t}\Phi
=-\delta\mathcal{W} + U
\quad\text{in }(0,T)\times S^2,
\end{align}
if for every $\varphi\in C^\infty_c((0,T)\times S^2,\R^3)$
there holds
\begin{align}
-\int_0^T\int_{S^2}
\Big\langle\Phi,\frac{\partial}{\partial t}
\varphi\Big\rangle\,
d\sigma_g\, dt
=-\int_0^T
\big(
\delta\mathcal{W}d\sigma_g, \varphi(t,\cdot)
\big)_{\mathcal{D}'}dt
+\int_0^T\int_{S^2}\big\langle U,\varphi\big\rangle\,
d\sigma_g\, dt.
\end{align}
\end{definition}

Our main result is the following.
\begin{theorem}
\label{thm:main}
There exists $\varepsilon_0>0$ so that
the Cauchy problem for the
conformal Willmore flow
\eqref{eq:Cauchy-conformal-Willmore}-\eqref{eq:conformal-flow-cond}
with initial datum in $\mathscr{W}^{\varepsilon_0}(S^2,\R^3)$
has a weak solution in
$\mathscr{W}^{\varepsilon_0,\delta}_{[0,T)}(S^2,\R^3)$
for some $\delta>0$,
assuming the initial datum in the sense of traces.
Such solution is smooth, exists for
all times and smoothly converges
to the standard embedding  $I$ of $S^2$ in $\R^3$.
Furthermore, if the initial datum is smooth,
such weak solution is also unique.
\end{theorem}

We can compare this result with the above metioned one
of \textsc{Kuwert} and \textsc{Sch\"atzle} \cite{MR1882663,MR1900754}.
They obtain,
in the smooth class, long-time existence, uniqueness and convergence to a round sphere for
the Cauchy problem of the normal flow
\eqref{eq:Cauchy-normal-Willmore}.
A central feature our result is that
the uniqueness of this smooth solution is in the broad class
of finite energy solutions,
and the fact that it converges exactly to the standard embedding.

We expect the solution to be unique
also if the initial datum is nonsmooth;
we plan to address this question in the future.
\medskip

The proof of the regularity part of Theorem \ref{thm:main}
shares evident similarities with
the corresponding one for the harmonic map flow
obtained in \cite{RiviereThesisIV}.
In that work, the core estimate that was obtained
for weak solutions of the harmonic map flow was of the form
\begin{align*}
\|u(t,\cdot)\|_{W^{2,2}}
\leq C\big(\|\partial_t u (t,\cdot)\|_{L^2} +1\big)
\quad\text{for a.e. } t,
\end{align*} 
which could then be squared and integrated in time
to yield higher regularity, and eventually smoothness
by the classical theory by \textsc{Struwe} 
\cite{MR826871,MR2431434}.
We shall obtain a similar result, namely an inequality
of the form
\begin{align*}
\|\Phi(t,\cdot)\|_{W^{4,2}}
\leq C\big(\|\e^{\lambda}\delta\mathcal{W}(t,\cdot)\|_{L^2} +1\big)
\quad\text{for a.e. } t,
\end{align*} 
for weak solutions of the conformal Willmore flow,
and likewise obtain higher regularity from it.
The overall procedure shall be however more technical. 

\subsubsection{Final Comments}
We have intentionally decided to work in a small-energy regime
in this first paper. 
We plan to consider more general scenarios
in future works, where more technical, localization/energy-concentration arguments will be
dealt with.

One has also to take into account
that, when the underlying surface is not a sphere, 
there is more than one conformal class,
so to properly work with a conformal Willmore flow,
one has to take into account the nontriviality of
the corresponding  Teichm\"uller space.
The
work of \textsc{Rupflin} and \textsc{Topping} 
\cite{MR3538152}
on the Teichm\"uller harmonic map flow also faces
the difficulty of ``following'' the conformal class along the flow.

In the future we plan to
determine whether the class
$\mathscr{W}^\varepsilon(S^2,\R^3)$
coincides or not with $\mathscr{E}(S^2,\R^3)$,
and whether the solution given by Theorem \ref{thm:main}
is unique also in the case of weak initial data.
We shall also seek to extend 
the argument  for branched weak initial data.

Finally, we plan to
carry an accurate study of singularities
(blow-up points, degeneration of conformal factor or conformal class...)
in forthcoming works.
For this we will likely build upon some of the work
already done on the subject,
\textsc{Mayer} and \textsc{Simonett} \cite{MR1877537},
\textsc{Blatt} \cite{MR2591055} and
\textsc{Chill}, \textsc{Fa\v{s}angov\'a} and \textsc{Sch\"atzle}
\cite{MR2495079} just to mention a few.

One of the questions relative to the parametric
approach of the Willmore flow is the following:
can the conformal class of a 
conformal Willmore flow 
 -- suitably normalized to remove any obvious gauge invariance --
degenerate in finite time?

\section{Preliminaries}
\label{sec:results-literature}
We state here known results that
shall play a key role in the paper.
\begin{theorem}
\label{thm:Wente-variant}
(\cite{MR1168525}, \cite{RiviereThesisIV}, \cite{MR1675268})
Let $\Omega\subseteq\R^2$ be a bounded domain,
let $a\in W^{1,2}(\Omega)$ and $f\in L^p(\Omega)$ for $1<p<2$.
Let $u\in W^{1,(2,\infty)}(\Omega)$%
\footnote{That is, $u\in W^{1,1}(\Omega)$
and $\nabla u$ is in the Lorentz space $L^{(2,\infty)}(\Omega)$.}
 be a solution to
\begin{align}
-\Delta u = \langle \nabla^\perp a, \nabla u\rangle + f
\quad\text{in }\Omega.
\end{align}
Then $u\in W^{2,p}_{\loc}(\Omega)$.
\end{theorem}

\begin{theorem}
[\cite{MR1913803}, \cite{MR3524220}]
\label{thm:eps-control-conformal-fact}
There exists $\varepsilon_0>0$ so that, if 
$\Phi\in \mathscr{E}(B_1,\R^n)$ is conformal with conformal factor $e^{\lambda}$
and so that
\begin{align}
\int_{B_1} |A|^2_g\,d\sigma_g
\leq \varepsilon_0,
\end{align}
then, if $C_{(2,\infty)}>0$ is constant so that
\begin{align}
\|d\lambda\|_{L^{(2,\infty)}(B_1)}\leq C_{(2,\infty)},
\end{align}
for any $0<r<1$ there holds
\begin{align}
\|d\lambda\|_{L^2(B_r)} + \|\lambda - \ell\|_{L^\infty(B_r)}
\leq C\int_{B_1}|A|^2_g\,d\sigma_g,
\end{align}
for some constants $\ell\in \R$ and
$C=C(r,C_{(2,\infty)})>0$.
\end{theorem}
\begin{remark}
By the triangle inequalily, without loss of generality 
we can take $\ell = \lambda(0)$  in the above estimate.
\end{remark}

\begin{theorem}[\cite{MR2169583,MR2232206}]
\label{thm:DeLellis-Muller}
Let $\mathcal{S}\subset\R^3$ be an immersed surface
with area 
$\mathcal{A}(\mathcal{S})=4\pi$
and let $g$ be the induced metric.
Then  there holds
\begin{align}
\label{eq:DeLellis-Muller-ineq}
\int_{\mathcal{S}}|A_{\text{sc}} -  g|^2_g\,d\sigma_g
\leq C\int_{\mathcal{S}}|A^\circ|^2_g\,d\sigma_g,
\end{align}
where 
$A_{\text{sc}}(\cdot\,,\cdot)=\langle A(\cdot\,,\cdot), N\rangle$ 
is the scalar second fundamental form of $\mathcal{S}$.
Furthermore, there is $\varepsilon_0>0$ so that if
\begin{align}
\mathcal{W}_0(\mathcal{S})
=\frac{1}{2}
\int_{\mathcal{S}}|A^\circ|^2_g\,d\sigma_{g}
\leq\varepsilon_0,
\end{align} 
there exists a conformal parametrization
$\Phi:S^2\to\R^3$ satisfying
\begin{align}
\label{eq:DeLellis-Muller-est}
\|\Phi - (c_{\mathcal{S}}+I)\|_{W^{2,2}(S^2)}+
\|\e^{\lambda}-1\|_{L^\infty(S^2)}
\leq C\sqrt{\mathcal{W}_0(\mathcal{S})},
\end{align}
for some vector $c_{\mathcal{S}}\in\R^3$ and for
an absolute constant $C>0$,
where
$\e^\lambda$ is the conformal factor of
the induced metric and 
$I$ is the standard immersion of $S^2$
into $\R^3$.
\end{theorem}
\begin{remark}
\label{rmk:minimality-avg}
From the minimality property
of the average:
\begin{align}
\Big\|f - \dashint_{S^2}f\,d\sigma\Big\|_{L^2(S^2)}
=\inf_{c\in\R}\|f-c\|_{L^2(S^2)},
\end{align}
and since $\int_{S^2} I\,d\sigma = 0$,
we may suppose
$c_{\mathcal{S}}=\dashint_{S^2} \Phi\,d\sigma_{S^2}$
in \eqref{eq:DeLellis-Muller-est}.
\end{remark}


\begin{theorem}
[\cite{MR2430975}, \cite{MR3524220}, \cite{MR3518329}]
\label{thm:Willmore-eqns}
The Willmore operator of an immersion
of a surface $\Phi:\Sigma\to\R^n$,
\begin{align}
\label{eq:Willmore-classical}
\delta\mathcal{W}
=\Delta^\perp_g H +  Q(A^\circ)H,
\end{align}
where $\Delta^\perp_g$ denotes the Laplace operator  on
the normal bundle of $\Phi(\Sigma)$
and 
\begin{align}
\label{eq:Willmore-Q}
 Q(A^\circ)H
=\langle A^\circ, \langle H,A^\circ\rangle\rangle_g
= g^{\mu\sigma}g^{\nu\tau}\langle
A^\circ_{\mu\nu}, \langle A^\circ_{\sigma\tau},H\rangle
\rangle,
\end{align}
may be written equivalently in divergence form as
\begin{align}
\label{eq:Willmore-divergence-form}
\delta\mathcal{W}
=\nabla^{*_g}\big( \nabla H 
+\langle A^\circ,H\rangle^{\sharp_g}
+\langle A,H\rangle^{\sharp_g}\big)
=\Delta_g H
+\nabla^{*_g}\big( 
\langle A^\circ,H\rangle^{\sharp_g}
+\langle A,H\rangle^{\sharp_g}\big).
\end{align}
When $n=3$, if
$w\in \Gamma(\Phi^\star(T\R^3)\otimes T^*\Sigma)$
denotes the vector-valued form along $\Phi$ given by
\begin{align}
\label{eq:w}
w= \nabla H 
+\langle A^\circ,H\rangle^{\sharp_g}
+\langle A,H\rangle^{\sharp_g}
=\nabla H
-2(\nabla H)^\top
-|H|^2d\Phi.
\end{align}
Then the following formulas hold true:
\begin{align}
\label{eq:Willmore-trasl}
\delta\mathcal{W} &= \nabla{}^{*_g}w,\\
\label{eq:Willmore-rot}
\Phi\times\delta\mathcal{W}
&=\nabla{}^{*_g}
\big(
-d\Phi\times H
+\Phi\times w
\big),\\
\label{eq:Willmore-dil}
\big\langle
\Phi,\delta\mathcal{W}
\big\rangle
&=d^{*_g}
\big\langle
\Phi, w
\big\rangle,\\
\label{eq:Willmore-inv}
\Phi\times(\Phi\times \delta \mathcal{W}) 
+\Phi\big\langle \Phi, \delta \mathcal{W}\big\rangle +4H&=
\nabla{}^{*_g}\big(
-2\Phi\times(d\Phi\times H)
+\Phi\times(\Phi\times w)
+\Phi\langle \Phi,w\rangle
\big).
\end{align}
\end{theorem}

\begin{remark}
\begin{enumerate}[label=(\roman*)]
\item 
Equations
\eqref{eq:Willmore-trasl}, 
\eqref{eq:Willmore-rot},
\eqref{eq:Willmore-dil}
and \eqref{eq:Willmore-inv}
follow from Noether's theorem
and the invariance of the Willmore
energy under translations, rotations, dilations
and inversions respectively.
To obtain \eqref{eq:Willmore-inv} 
(since inversions do not form \textit{per se} a
1--parameter family of transformations),
one applies Noether's theorem to the local family
$ F(s,y) = \mathcal{I} \, \circ \, \tau_{sv} \circ \mathcal{I}(y)$,
where  $\mathcal{I}(y) = \frac{y}{|y|^2}$
denotes the inversion with respect to the unit sphere
and $\tau_v(y) = y + v$ is the translation
of vector $v\in\R^3$.
This transformation is generated by the vector field
$\mathcal{X}(y) = v|y|^2 -2 y\langle y, v\rangle$
and gives rise to the law \eqref{eq:Willmore-inv}.
\item The equivalence between the first 
and second expressions for $w$ in \eqref{eq:w}
is obtained by noting that,
since $H$ is a normal vector, there holds
\begin{align}
(\nabla H)^\top
& = g^{\mu\nu}
\langle
\partial_\alpha H , \partial_\mu \Phi
\rangle \partial_\nu\Phi \otimes dx^\alpha \\
& = - g^{\mu\nu} \langle
H , \partial^2_{\alpha\mu} \Phi
\rangle \partial_\nu\Phi \otimes dx^\alpha \\
& = - g^{\mu\nu} \langle
H , A_{\alpha\mu}
\rangle \partial_\nu\Phi \otimes dx^\alpha \\
& = -\langle H, A \rangle ^{\sharp_g},
\end{align} 
and since one similarly has
$|H|^2 d\Phi = \langle H, H g\rangle^{\sharp_g}$,
there holds
\begin{align}
2(\nabla H)^\top
+|H|^2 d\Phi
& = -2\langle H, A \rangle^{\sharp_g}
+ \langle H, H g\rangle^{\sharp_g}
= -\langle H, A \rangle^{\sharp_g}
 - \langle H, A^\circ \rangle^{\sharp_g}.
\end{align}
\end{enumerate}
\end{remark}

\begin{theorem}
[\cite{MR1882663,MR1900754,10.1093/imrn/rnaa015}]
\label{thm:area-est}
There exists an $\varepsilon_0=\varepsilon_0(n)>0$ 
so that, if ${\Phi:[0,T)\times S^2\to \R^n}$
is a smooth normal Willmore flow
\begin{align}
\frac{\partial}{\partial t}\Phi = -\delta\mathcal{W},
\end{align}
with smooth initial datum $\Phi(0,\cdot)=\Phi_0$ and
$\mathcal{W}_0(\Phi_0)\leq\varepsilon_0$,
then:
\begin{enumerate}[label=(\roman*)]
\item  its area satisfies
\begin{align}
\label{eq:area-est}
|\mathcal{A}(\Phi(t,\cdot)) - \mathcal{A}(\Phi_0)|
\leq C\mathcal{A}_0(\Phi_0)\mathcal{W}_0(\mathcal{S}),
\end{align}
for a constant $C=C(n)>0$.
\item Its barycenter
$\mathcal{C}(\Phi)=\dashint_{S^2}\Phi\, d\sigma_{g}$
satisfies
\begin{align}
|\mathcal{C}(\Phi(t,\cdot))-\mathcal{C}(\Phi_0)|\leq C\mathcal{W}_0(\mathcal{S})
\end{align}
for a constant $C=C(n)>0$.
\item $\Phi$ exists for all times and smoothly converges
to a round sphere $t\to\infty$.
\end{enumerate}
\end{theorem}

\section{Consequences of DeLellis-M\"uller's Theorem}
\label{sec:DLM-consequences}
%

The proof of Proposition \ref{prop:well-balanced-fam}
will follow from two lemmas.

\begin{lemma}
\label{lemma:key-IFT}
The function $F:\Aut(S^2)\to \R^6$ given by
\begin{align}
\label{eq:key-IFT-F}
F(\psi) =
(F_1(\psi),F_2(\psi))=
\Big(
\int_{S^2} I\, \psi^\star d\sigma,
\int_{S^2} (I\circ\,\psi)\times I\,d\sigma
\Big),
\end{align}
where $\psi^\star d\sigma$ denotes the pullback
of area element $d\sigma$ via $\psi$,
is differentiable and $dF(e)$ is an isomorphism. 
\end{lemma}

\begin{bfproof}
Differentiability  follows
since $F$ is composition of smooth functions and operations.
We shall now use the language of differential
forms and so along this proof it is convenient to
temporarily change our notation for the area element
from $d\sigma$ to $\omega_{S^2}$.
\medskip

Recall that $\omega_{S^2}=\omega_{\R^3}\llcorner N
=\omega_{\R^3}(N,\cdot,\cdot)$,
where $N$ Gauss map of $S^2$.
More explicitly, since 
$N(y) = y$,
we have the formula
\begin{align}
\omega_{S^2}
=(dy^1\wedge dy^2\wedge dy^3)(y,\cdot,\cdot)
=y^1dy^2\wedge dy^3
-y^2dy^1\wedge dy^3
+y^3dy^1\wedge dy^2.
\end{align}
Now $\Aut(S^2)$ has dimension 6 as a manifold
and we consider the basis for the tangent space $T_e\Aut(S^2)$
given by the vector
fields generating, respectively, rotations
and ``spherical dilations'' about the coordinate axes:
\begin{align}
Z_1(y) &=
(0,-y^3,y^2), &
Z_2(y) &=
(y^3,0,-y^1),&
Z_3(y) &=
(-y^2,y^1,0),\\[0.5em]
Z_4(y) &=
e_1-y^1 y, &
Z_5(y) &=
e_2-y^2 y,&
Z_6(y) &=
e_3-y^3 y.&
\end{align}
We shall prove  that 
\begin{align}
\label{eq:key-IFT-matrix}
\Big(\frac{\partial F^j}{\partial Z_a}(e)\Big)_{1\leq a,j\leq 6}
=-\frac{8\pi}{3}\mathbf{1}_{6\times 6}.
\end{align}
To compute $\partial_X F(e)$, if $\Phi^X$ denotes the local flow
of the vector field $X$, we have to evaluate
\begin{align}
\frac{\partial}{\partial X}F(e)
=\frac{d}{d t}F(\Phi^X(t,\cdot))\bigg|_{t=0}.
\end{align}
Let us look at $F_1$.
By Cartan's formula, since $d\omega_{S^2}=0$, it is
\begin{align}
\frac{\partial}{\partial t}\Phi^{X}(t,\cdot) \omega_{S^2}
\Big|_{t=0} = \mathcal{L}_{X}\omega_{S^2}
= d(X\llcorner\omega_{S^2}),
\end{align}
where $\mathcal{L}_X$ denotes the Lie derivative with respect
to $X$. Since $Z_a$'s for $a=1,2,3$ generate isometries,
$\mathcal{L}_{Z_a}\omega_{S^2}=0$
and hence \emph{a fortiori}
\begin{align}
\frac{d}{dt}F_1\big(\Phi^{Z_a}(t,\cdot)\big)\bigg|_{t=0}
=\int_{S^2} I\,\mathcal{L}_{Z_a}\omega_{S^2}
=0
\quad\quad\text{for } a=1,2,3.
\end{align}
As for the $Z_a$'s for $a=4,5,6$, we see that
{\small
\begin{align}
Z_4\llcorner\omega_{S^2}
&=-y^2dy^3+y^3dy^2,&
Z_5\llcorner\omega_{S^2}
&=y^1dy^3-y^3dy^1,&
Z_6\llcorner\omega_{S^2}
&=-y^1dy^2+y^2dy^1,\\[0.5em]
d(Z_4\llcorner\omega_{S^2})
&=-2dy^2\wedge dy^3,&
d(Z_5\llcorner\omega_{S^2})
&=2dy^1\wedge dy^3,&
d(Z_6\llcorner\omega_{S^2})
&=-2 dy^1\wedge d y^2,
\end{align}
}
and hence with Stokes' theorem we get
{\small
\begin{align}
\frac{d}{dt}F_1\big(\Phi^{Z_4}(t,\cdot)\big)\bigg|_{t=0}
&=\int_{S^2} I\mathcal{L}_{Z_4}\omega_{S^2}\\
&=\Big(
\int_{S^2} y^1(-2dy^2\wedge dy^3),
\int_{S^2}y^2 (-2dy^1\wedge dy^3),
\int_{S^2}y^3(-2 dy^1\wedge d y^2)
\Big)\\
&=-\frac{8\pi}{3}(1,0,0),
\end{align}
}
and similarly
\begin{align}
\frac{d}{dt}F_1\big(\Phi^{Z_5}(t,\cdot)\big)\bigg|_{t=0}
=-\frac{8\pi}{3}(0,1,0),\qquad
\frac{d}{dt}F_1\big(\Phi^{Z_6}(t,\cdot)\big)\bigg|_{t=0}
=-\frac{8\pi}{3}(0,0,1).
\end{align}
Now we consider $F_2$.
Since we can write
\begin{align}
F_2\big(\Phi^{X}(t,\cdot)\big)
=\int_{S^2}
\Phi^{X}(t,\cdot)\times I\,d\sigma,
\end{align}
it is
\begin{align}
\frac{d}{dt}F_2\big(\Phi^{X}(t,\cdot)\big)\bigg|_{t=0}
=\int_{S^2}
X\times I\,d\sigma,
\end{align}
and thus one directly computes that
for $a=1,2,3$ it is
\begin{align}
\int_{S^2}
Z_a\times I\,d\sigma
&=-\frac{8\pi}{3}e_a,
\qquad
\int_{S^2}
Z_{a+3}\times I\,d\sigma
=0.
\end{align}
Putting together all these computations yields 
\eqref{eq:key-IFT-matrix} and hence the thesis.
\end{bfproof}

\begin{lemma}
\label{lemma:WB-exists}
There exists $E>0$ with the following property.
For any $\eta_1\leq E$
there exist $\eta_2>0$ so that,
if $\mathcal{S}\subset\R^3$ is an immersed surface
with area $\mathcal{A}(\mathcal{S})=4\pi$
and $\Phi:S^2\to\mathcal{S}$ is a conformal immersion
with conformal factor $\e^{\lambda}$ and
$c\in\R^3$ is a vector so that
\begin{align}
\label{eq:WB-a}
\|\Phi - I-c\|_{W^{2,2}(S^2)}
+\|e^\lambda -1\|_{L^\infty(S^2)}
\leq \eta_1,
\end{align}
then there exists a conformal self-map 
$\psi\in \Aut(S^2)$ so that
$\Psi = \Phi\circ \psi$ is well-balanced and
$\psi$ is the unique self-map with such property in the 
Riemannian ball $B_{\eta_2}(e)\subset \Aut(S^2)$.
In addition,
if $\e^{\nu}$ denotes the conformal factor of $\Psi$,
there holds
\begin{align}
\label{eq:WB-b}
\|\Psi - I -c\|_{W^{2,2}(S^2)}
+\|\e^{\nu}-1\|_{L^\infty(S^2)}
\leq 2\eta_1.
\end{align}
\end{lemma}

\begin{bfproof}
Let $\mathcal{F}:W^{2,2}(S^2,\R^3)\times \Aut(S^2)\to \R^6$
be given by
\begin{align}
\quad
\mathcal{F}(f,\psi)
=\Big(
\int_{S^2}I\,\frac{1}{2}|d(f\circ\psi)|^2d\sigma,
\int_{S^2}(f\circ\psi)\times I\,d\sigma
\Big).
\end{align}
Note that this definition makes sense for every $f\in W^{2,2}(S^2)$
and, if $\Phi$ is a conformal immersion
$\mathcal{F}(\Phi,e)=0$ means that $\Phi$ is well-balanced as
in Definition \ref{def:well-balanced}.
Moreover  $\mathcal{F}(I,\cdot)=0$
coincides with $F$ given in \eqref{eq:key-IFT-F}.
Finally $\mathcal{F}$ is invariant by translations
in its first component: $\mathcal{F}(f,\cdot) = \mathcal{F}(f+k,\cdot)$
for every $k\in\R^3$.
\medskip

As a consequence of Lemma \ref{lemma:key-IFT}, 
$d_\psi \mathcal{F}(I,e)=d_\psi \mathcal{F}(I,\cdot)(e)$
is an isomorphism,
and hence
by the implicit function theorem,
there exists $E,\eta_2>0$ so that
if \eqref{eq:WB-a} holds for $\eta_1\leq E$,
there is a unique $\psi = \psi_\Phi$ in the
Riemannian ball $B_{\eta_2}(e)\subset \Aut(S^2)$
 so that $\mathcal{F}(\Phi-c,\psi)=\mathcal{F}(\Phi,\psi)=0$,
i.e. so that $\Psi=\Phi\circ\psi$ is well-balanced
(recall also Remark \ref{rmk:balance-transl}).
\medskip

Finally, since $\psi$ is biholomorphic,
$\forall N \in \N$ we can estimate
${\sum_{k=1}^N\dist(\nabla^k\psi,\nabla^k e)
\leq C_N\dist (\psi,e)}$
for some $C_N>0$ independent of $\psi$.
So
\eqref{eq:WB-a} holds,  by the triangle inequality
and the continuity of the Lebesgue integral,
\begin{align}
&\phantom{{}={}}
\|\Psi - I-c\|_{W^{2,2}(S^2)}
+\|\e^{\nu}-1\|_{L^\infty(S^2)}\\
& = 
\|\Phi\circ\psi-I-c\|_{W^{2,2}(S^2)}
+\|\tfrac{1}{\sqrt{2}}|d\psi|\e^{\lambda\circ\psi}-1\|_{L^\infty(S^2)}\\
& = 
\|(\Phi- c)\circ\psi -I\|_{W^{2,2}(S^2)}
+\|\tfrac{1}{\sqrt{2}}|d\psi|\e^{\lambda\circ\psi}-1\|_{L^\infty(S^2)}\\
&=\eta_1 + o(1)
\quad\text{as } 
\dist(\psi,e)\to 0,
\end{align}
and, since $\|\Psi-c\|_{L^2(S^2)}$
and $\|e^{\lambda}\|_{L^\infty(S^2)}$ are
uniformly bounded, the remainder $o(1)$
can be taken uniform in $\Phi,\Psi,\psi,c$
and hence, choosing $\eta_2$ sufficiently small
we obtain to \eqref{eq:WB-b}.
\end{bfproof}

\begin{bfproof}[Proof of Proposition \ref{prop:well-balanced-fam}]
It suffices to prove the thesis
for $\mathcal{W}_0(\mathcal{S})\leq\varepsilon_0$
sufficiently small.
For part (i),
combine
Theorem \ref{thm:DeLellis-Muller}
and Lemma \ref{lemma:WB-exists}.

For part (ii),
Let 
$\varepsilon_0>0$
be sufficiently small so that 
\begin{align}
C\varepsilon_0
\leq \frac{1}{2}E,
\end{align}
where $E$ is as in Lemma \ref{lemma:WB-exists}
and $C$ is the constant of Theorem \ref{thm:DeLellis-Muller},
and let 
$\Phi:S^2\to\mathcal{S}$ be the conformal
parametrization given by that theorem.
\medskip
By Lemma \ref{lemma:WB-exists}
there exists
a unique choice of  $\alpha=\alpha_{\Phi}$ in $B_{\eta_2}(e)\subset\Aut(S^2)$
with $\mathcal{F}(\Phi, \alpha)=0$ i.e. $\Phi'=\Phi\circ\alpha$ is well-balanced
and
\begin{align}
\|\Phi' - I-c\|_{W^{2,2}(S^2)}
+\|e^{\lambda'} -1\|_{L^\infty(S^2)}
\leq 
2C\sqrt{\mathcal{W}_0(\Phi')}
\leq \sqrt{E}.
\end{align}
So now if $\delta$ is taken so that
\begin{align}
\delta\leq\frac{1}{2}\sqrt{E},
\end{align}
since $\Psi$ is already well-balanced, by uniqueness it must be
$\Phi'=\Psi$, and
the thesis follows also for the local uniqueness part, with 
$\mathcal{O}=B_{\eta_2}(e)$.
\end{bfproof}

The following simple consequence
of Theorem \ref{thm:DeLellis-Muller}
will also be needed later.
\begin{lemma}
\label{lemma:DeLellis-Muller-local}
If $\Phi:S^2\to \R^3$ is a conformal immersion
with conformal factor $e^{\lambda}$ and
$B_r(x_0)\subset S^2$ is a disk of radius $r$ 
(in the standard metric of $S^2$),
there holds  
\begin{align}
\label{eq:DeLellis-Muller-local}
\int_{B_r(x_0)}
|A|^2_g\,d\sigma_g
&\leq C\Big(\int_{S^2}
|A^\circ|^2_g\,d\sigma_g+
\e^{4C_0}r^2\Big),
\end{align}
where $C_0 = \|\lambda\|_{L^\infty(S^2)}$
and $C>0$ is an absolute constant.
\end{lemma}

\begin{bfproof}
It is a consequence of 
\eqref{eq:DeLellis-Muller-ineq} applied to
the immersed surface $\mathcal{S} = a\,\Phi(S^2)$, where
 $a=\sqrt{\frac{4\pi}{\mathcal{A}(\Phi)}}$ and
$\mathcal{A}(\Phi)=\int_{S^2}\e^{2\lambda}\,d\sigma$ is the area of $\Phi(S^2)$.
Indeed, since
\begin{align}
4\pi\e^{-2C_0}\leq \mathcal{A}(\Phi)\leq 4\pi\e^{2C_0},
\end{align}
it follows that
$e^{-C_0}\leq a \leq \e^{C_0}$
and we can estimate
\begin{align}
\int_{B_r(x_0)}|a\, g|^2_gd\sigma_g 
=a^2 \int_{B_r(x_0)}2\e^{2\lambda}d\sigma
\leq C\e^{4C_0}r^2.
\end{align}
and thus
\begin{align}
\int_{B_r(x_0)}
|A|^2_g\,d\sigma_g
&\leq
2\int_{B_r(x_0)}
|A_{\text{sc}}-a g|^2_g\,d\sigma_g
+2\int_{B_r(x_0)}
|a\, g|^2_g\,d\sigma_g\\
&\leq
C \int_{B_r(x_0)}
|A^\circ|^2_g\,d\sigma_g
+C\e^{4C_0}r^2,
\end{align}
which proves \eqref{eq:DeLellis-Muller-local}.
\end{bfproof}

\section{Estimates of Elliptic Type}
\label{sec:ellptic}
This section is dedicated to prove the following.
\begin{theorem}
\label{thm:quant-L2-reg-Willmore}
Let $\Phi\in \mathscr{E}(B_1,\R^3)$ be conformal
with conformal factor $\e^\lambda$ and
Willmore operator $\delta\mathcal{W}$ in $L^2(B_1)$.
Then $\Phi\in W^{4,2}_{\loc}(B_{1})$, and furthermore
if $C_{(2,\infty)}>0$ is constant so that
\begin{align}
\|d\lambda\|_{L^{(2,\infty)}(B_1)}\leq C_{(2,\infty)},
\end{align}
there exists an $\varepsilon_0>0$
depending only on $C_{(2,\infty)}$
so that if
\begin{align}
\int_{B_1}|A|^2_g\,d\sigma_g\leq\varepsilon_0,
\end{align}
then
the following estimate holds:
\begin{align}
\label{eq:quant-L2-est-Willmore}
\|d\Phi\|_{W^{3,2}(B_{1/2})}
\leq C\big(\|\e^{4\lambda}\delta\mathcal{W}\|_{L^2(B_1)}
+\|\e^{\lambda}\|_{L^2(B_1)}\big),
\end{align}
where $C=C(C_{(2,\infty)})>0$.
\end{theorem}

\begin{remark}
The fact that the estimate \eqref{eq:quant-L2-est-Willmore}
does not include $\|\Phi\|_{L^2(B_{1/2})}$
on the left hand-side
is motivated by the translation invariance 
of all the quantities on the right-hand side. 
\end{remark}

\subsection{Qualitative Estimates}
\begin{proposition}
\label{prop:bootstrap-cond-Willmore}
Let $1<q<\infty$ and let 
$\Phi\in\mathscr{E}(B_1,\R^n)$ be a conformal with
Willmore operator $\delta\mathcal{W}$ in $L^q(B_1)$.
It suffices to know that $H\in L^p(B_1)$ for some $p>2$
to deduce that $\Phi\in W^{4,q}_{\loc}(B_1)$.
\end{proposition}

\begin{bfproof}
We may certainly assume $2<p<4$.
Since
\begin{align}
\Delta\Phi = 2\e^{2\lambda}H
\in L^{p},
\end{align}
elliptic regularity theory gives that 
$\Phi\in W^{2,p}_{\loc}$,
and this in turn implies
\begin{align}
A = (\nabla^2\Phi)^\perp \in L^{p}_{\loc}.
\end{align}
Looking at \eqref{eq:Willmore-divergence-form}:
\begin{align}
-\Delta H 
= \nabla^*
(\langle A^\circ,H\rangle^{\sharp_g} +\langle A,H\rangle^{\sharp_g})
-\e^{2\lambda}\delta\mathcal{W},
\end{align}
we see that
\begin{align}
\langle A^\circ,H\rangle^{\sharp_g} +\langle A,H\rangle^{\sharp_g}
=\e^{-2\lambda}\big(\langle A^\circ,H\rangle^{\sharp} +\langle A,H\rangle^{\sharp}\big)
\in L^{\frac{p}{2}}_{\loc},
\end{align}
and, since $q>1$, 
we have in particular that $\delta\mathcal{W}\in W^{-1,2}$;
thus $\Delta H\in W^{-1,\frac{p}{2}}_{\loc}$,
whence elliptic regularity and Sobolev embedding
give
\begin{align}
H\in W^{1,\frac{p}{2}}_{\loc}
\hookrightarrow L^{\left(\frac{p}{2}\right)^*}_{\loc},
\end{align}
and from this, it follows that
$
A\in L^{\left(\frac{p}{2}\right)^*}_{\loc}.
$
Then
\begin{align}
\langle A^\circ,H\rangle^{\sharp_g} +\langle A,H\rangle^{\sharp_g}
\in L^{\frac{1}{2}\left(\frac{p}{2}\right)^*}_{\loc},
\end{align}
whence $\Delta H\in W^{-1,\frac{1}{2}\left(\frac{p}{2}\right)^*}_{\loc}$,
so by elliptic regularity
\begin{align}
H\in W^{1,\frac{1}{2}\left(\frac{p}{2}\right)^*}_{\loc}.	
\end{align}
This process can be iterated, and
since the sequence
$p$, $\big(\frac{p}{2}\big)^*$, 
$\big(\frac{1}{2}\big(\frac{p}{2}\big)^*\big)^*$,
$\big(\frac{1}{2}\big(\frac{1}{2}\big(\frac{p}{2}\big)^*\big)^*\big)^*$,$\ldots$
is strictly monotone increasing and unbounded,
after a finite number (depending on $p$) of steps
we get that
\begin{align}
\Big(\frac{1}{2}\Big(\frac{1}{2}\Big(\cdots
\Big(\frac{p}{2}\Big)^*\cdots\Big)^*\Big)^*\Big)^*
\geq 2.
\end{align}
We then deduce that $-\Delta H\in W^{-1,2}_{\loc}$,
and thus that $H\in W^{1,2}_{\loc}$.
By Sobolev embedding, 
this yields $H\in L^r_{\loc}$ for every $r<\infty$,
and in turn elliptic estimates give
\begin{align}
\Phi\in W^{2,r}_{\loc}
\quad\forall\,r<\infty,
\end{align}
hence also $A\in L^r_{\loc}$ for every $r<\infty$.
From Liouville equation
\begin{align}
-\Delta \lambda = \e^{2\lambda}K,
\end{align}
since $|K|\leq C|A|^2$, we have $\Delta\lambda\in L^r_{\loc}$
and hence $\lambda\in W^{2,r}_{\loc}$ for every $r<\infty$.
With this we infer that in fact
$A, H\in W^{1,r}_{\loc}$ and so
\begin{align}
\langle A^\circ,H\rangle^{\sharp_g} +\langle A,H\rangle^{\sharp_g}
\in W^{1,r}_{\loc},
\end{align}
Thus
$-\Delta H\in L^q_{\loc}$, which implies
$
H\in W^{2,q}_{\loc}
$
and hence
(again since $\lambda\in W^{2,r}_{\loc}$)
that $\Phi\in W^{4,q}_{\loc}$.
\end{bfproof}

\begin{lemma}
\label{lemma:w-scal-dPhi}
With the same hypotesis as in Theorem \ref{thm:Willmore-eqns},
the vector-valued form \eqref{eq:w}
satisfies $\langle w,d\Phi\rangle_g=0$.
\end{lemma}
\begin{bfproof}
Since $|d\Phi|_g^2=2$ and $H$ is a normal vector field,
we have
\begin{align}
\left\langle
w,d\Phi
\right\rangle_g
&=\left\langle
\nabla H,d\Phi
\right\rangle_g
-2\left\langle
(\nabla H)^\top,d\Phi
\right\rangle_g
-|H|^2|d\Phi|_g^2\\
&=-\left\langle
\nabla H,d\Phi
\right\rangle_g
-2|H|^2\\
&=-\frac{1}{2}\left\langle
\nabla \Delta_g\Phi,d\Phi
\right\rangle_g
-2|H|^2\\
&=\frac{1}{2}\left\langle
\Delta_g\Phi,\nabla^* d\Phi
\right\rangle_g
-2|H|^2\\
&=\frac{1}{2}|\Delta_g\Phi|^2 - \frac{1}{2}|\Delta_g\Phi|^2
=0.
\qedhere
\end{align}
\end{bfproof}

\begin{proposition}
\label{prop:Willmore-inhom-system}
$\Phi:B_1\to\R^3$ be an immersion.
Consider the Hodge decomposition
\begin{align}
\label{eq:Hodge-w}
w=d\mathscr{L} + *dL,
\end{align}
and consider further the Hodge decompositions
\begin{align}
\label{eq:Hodge-R}
-d\Phi\times H - (*d\Phi)\times L &= d\mathscr{R} + *dR,\\
\label{eq:Hodge-S}
-\langle *d\Phi,L\rangle &= d\mathscr{S} + *dS.
\end{align}
Then the following relations hold:
\begin{align}
\label{eq:Willmore-inh-system-1}
\Delta_g\Phi
&=d\Phi\times
\big(d\mathscr{R} + *dR\big)
+\big\langle d\Phi,d\mathscr{S} + *dS\big\rangle\\
\label{eq:Willmore-inh-system-2}
\Delta_g R
&=dN\times(d\mathscr{R} + *dR)
-\langle dN, d\mathscr{S} + *dS\rangle
+(*d\Phi)\times d\mathscr{L},\\
\label{eq:Willmore-inh-system-3}
\Delta_g S
&=\langle dN, d\mathscr{R} + *dR\rangle 
+ \langle *d\Phi, d\mathscr{L}\rangle,
\end{align}
where $N$ is the Gauss map of $\Phi$.
\end{proposition}

\begin{bfproof}
Note first that, if the Hodge decompositions 
\eqref{eq:Hodge-w}, \eqref{eq:Hodge-R}, \eqref{eq:Hodge-S}
hold, we necessarily have
\begin{align}
\label{eq:mathscr-L}
\Delta_g\mathscr{L}&=d^*w=\delta \mathcal{W},\\
\label{eq:mathscr-R}
\Delta_g \mathscr{R} &=-d\Phi\times d\mathscr{L},\\
\label{eq:mathscr-S}
\Delta_g \mathscr{S} &=- \langle d\Phi,d\mathscr{L}\rangle.
\end{align}
With \eqref{eq:Willmore-inv} and the identity%
\footnote{
here, the notation is:
$
d\Phi\times dH
=g^{\mu\nu}\partial_{\mu}\Phi \times \partial_\nu H
$,
so that 
$
\Phi\times d\Phi\times dH
=g^{\mu\nu}(\Phi\times\partial_{\mu}\Phi \times \partial_\nu H).
$
Similarly for
$w\times dH$ and $\Phi\times d\Phi\times w$.
}
\begin{align}
d^*(-\Phi\times d\Phi\times H)
=2H - \Phi\times d\Phi\times dH
=2H - \Phi\times d\Phi\times w,
\end{align}
we see that there holds
\begin{align}
&\phantom{{}={}}\Phi\times(\Phi\times \delta\mathcal{W}) 
+ \Phi\langle \Phi,\delta\mathcal{W}\rangle
+4H\\
&=d^*\big(
-2\Phi\times(d\Phi\times H)
+\Phi\times(\Phi\times w) + \Phi\langle\Phi,w\rangle
\big)\\
&=d^*\big(-\Phi\times(d\Phi\times H)\big)\\
&\phantom{{}={}}
+d^*\Big(
\Phi\times\big(
-d\Phi\times H+\Phi\times (d\mathscr{L} + *dL)
\big) 
+ \Phi\big\langle\Phi,d\mathscr{L}+*dL\big\rangle
\Big)\\
&=2H -\Phi\times(d\Phi\times w)\\
&\phantom{{}={}}
+d^*\Big(
\Phi\times\big(
-d\Phi\times H+\Phi\times d\mathscr{L} +*d(\Phi\times L) 
-(*d\Phi)\times L\big)
\Big)\\
&\phantom{{}={}}
+ d^*\Big(
\Phi\big(
\langle\Phi,d\mathscr{L}\rangle
+*d\langle \Phi,L\rangle
-\langle*d\Phi,L\rangle
\big)
\Big)\\
&=2H -\Phi\times(d\Phi\times w)\\
&\phantom{{}={}}
+d^*\Big(
\Phi\times\big(
\Phi\times d\mathscr{L} 
+d\mathscr{R} + *dR +*d(\Phi\times L) 
\big)
\Big)\\
&\phantom{{}={}}
+ d^*\Big(
\Phi\big(
\langle\Phi,d\mathscr{L}\rangle
+d\mathscr{S} + *dS
+*d\langle \Phi,L\rangle
\big)
\Big).
\end{align}
Now, on the one hand we have,
from \eqref{eq:mathscr-R},
\begin{align}
&\phantom{{}={}}
d^*\Big(
\Phi\times\big(
\Phi\times d\mathscr{L} 
+d\mathscr{R} + *dR +*d(\Phi\times L) 
\big)
\Big)\\
&=d\Phi\times
\big(
\Phi\times d\mathscr{L} 
+d\mathscr{R} + *dR +*d(\Phi\times L) 
\big)+\Phi\times
\big(
d\Phi\times d\mathscr{L}
+\Phi\times \Delta_g\mathscr{L} 
+\Delta_g\mathscr{R}
\big)\\
&=d\Phi\times
\big(
\Phi\times d\mathscr{L} 
+d\mathscr{R} + *dR +*d(\Phi\times L) 
\big)+\Phi\times
(\Phi\times \delta \mathcal{W})
\end{align}
on the other hand from \eqref{eq:mathscr-S} it follows that
\begin{align}
&\phantom{{}={}}
d^*\Big(
\Phi\big(
\langle\Phi,d\mathscr{L}\rangle
+d\mathscr{S} + *dS
+*d\langle \Phi,L\rangle
\big)
\Big)\\
&=\big\langle
d\Phi,
\langle\Phi,d\mathscr{L}\rangle
+d\mathscr{S} + *dS
+*d\langle \Phi,L\rangle
\big\rangle
+\Phi\big(
\langle d\Phi,d\mathscr{L}\rangle
+\langle\Phi, \Delta_g\mathscr{L}\rangle
+\Delta_g\mathscr{S}
\big)\\
&=\big\langle
d\Phi,
\langle\Phi,d\mathscr{L}\rangle
+d\mathscr{S} + *dS
+*d\langle \Phi,L\rangle
\big\rangle
+\Phi\langle\Phi,\delta\mathcal{W}\rangle
\end{align}
thus we deduce
\begin{align}
2H&=-\Phi\times(d\Phi\times w)\\
&\phantom{{}={}}
+d\Phi\times
\big(
\Phi\times d\mathscr{L} 
+d\mathscr{R} + *dR +*d(\Phi\times L) 
\big)\\
&\phantom{{}={}}+\big\langle
d\Phi,
\langle\Phi,d\mathscr{L}\rangle
+d\mathscr{S} + *dS
+*d\langle \Phi,L\rangle
\big\rangle\\
&=-\Phi\times(d\Phi\times w)
+d\Phi\times
\big(d\mathscr{R} + *dR\big)
+\big\langle d\Phi,d\mathscr{S} + *dS\big\rangle\\
&\phantom{{}={}}
+d\Phi\times
\big(
\Phi\times d\mathscr{L}
+*d(\Phi\times L)
\big)
+ \big\langle
d\Phi,
\langle\Phi,d\mathscr{L}\rangle
+*\langle \Phi,L\rangle\big\rangle
\end{align}
By definition of $\mathscr{L}$ and $L$,
with Lemma \ref{lemma:w-scal-dPhi} the last line in the above
expression is 
\begin{align}
&\phantom{{}={}}
d\Phi\times
\big(
\Phi\times d\mathscr{L}
+*d(\Phi\times L)
\big)
+ \big\langle
d\Phi,
\langle\Phi,d\mathscr{L}\rangle
+*\langle \Phi,L\rangle\big\rangle\\
&=d\Phi\times
\big(
\Phi\times (d\mathscr{L} + *dL)
+(*d\Phi)\times L
\big)
+ \big\langle
d\Phi,
\langle\Phi,d\mathscr{L} + *d L\rangle
+\langle *d\Phi,L\rangle\big\rangle\\
&=d\Phi\times
\big(
\Phi\times w
+(*d\Phi)\times L
\big)
+ \big\langle
d\Phi,
\langle\Phi,w\rangle
+\langle *d\Phi,L\rangle\big\rangle\\
&=d\Phi\times
\big( \Phi\times w \big)
+ \big\langle
d\Phi,
\langle\Phi,w\rangle\big\rangle
+d\Phi\times\big((*d\Phi)\times L)
+\big\langle d\Phi, \langle *d\Phi, L\rangle\big\rangle\\
&=\Phi\times(d\Phi\times w),
\end{align}
and this yields \eqref{eq:Willmore-inh-system-1}.
Next, since $\langle d\Phi\times H,N\rangle=0$, we see that
\begin{align}
\big\langle d\mathscr{R} + *dR,N\big\rangle
&=\big\langle -d\Phi\times H -(*d\Phi)\times L, N\big\rangle\\
&=-\big\langle N\times (*d\Phi),L\big\rangle\\
&=-\big\langle d\Phi,L\big\rangle\\
&=-*d\mathscr{S} +dS
\end{align}
and similarly, using the rules of the vector product,
we have
\begin{align}
\big(d\mathscr{R} + *dR\big)\times N
&=(-d\Phi\times H - (*d\Phi)\times L\big)\times N\\
&=N\times\big(d\Phi\times H\big) + N\times\big((*d\Phi)\times L\big)\\
&=-H\times(N\times d\Phi)
-d\Phi\times(H\times N)\\
&\phantom{{}=}-L\times\big(N\times(*d\Phi)\big)
-(*d\Phi)\times(L\times N)\\
&=H\times(*d\Phi) - L\times d\Phi
-(*d\Phi)\times (L\times N)\\
&=-(*d\Phi)\times H + d\Phi\times L + N\langle *d\Phi,L\rangle\\
&=*d\mathscr{R} - dR - N(d\mathscr{S} + *dS),
\end{align}
and thus, codifferentiating these identities we have
\begin{align}
\Delta_g\mathscr{R}\times N
+\big(d\mathscr{R} + *dR\big)\times dN
&=-\Delta_g R
-\big\langle dN, d\mathscr{S} + *dS\big\rangle
-N\Delta_g \mathscr{S},\\
\big\langle \Delta_g \mathscr{R}, N\big\rangle
+\big\langle d\mathscr{R} + *dR,dN\big\rangle
&=\Delta_g S.
\end{align}
It now suffices to notice that, by definition of 
$\mathscr{R}$ and $\mathscr{S}$ it is
\begin{align}
N\Delta_g\mathscr{S} &= -N\langle d\Phi, d\mathscr{L}\rangle\\
\Delta_g \mathscr{R}\times N
&=-(d\Phi\times d\mathscr{L})\times N
=N\times(d\Phi\times d\mathscr{L})
=\big\langle d\Phi,\langle N,d\mathscr{L}\rangle\big\rangle,\\
\big\langle \Delta_g\mathscr{R}, N\big\rangle
&=-\langle d\Phi\times d\mathscr{L},N\rangle
=-\langle N\times d\Phi, d\mathscr{L}\rangle
=\langle *d\Phi,d\mathscr{L}\rangle
\end{align}
and in particular
\begin{align}
\Delta_g\mathscr{R}\times N + N\Delta_g \mathscr{S}
&=\big\langle d\Phi,\langle N,d\mathscr{L}\rangle\big\rangle
-N\langle d\Phi, d\mathscr{L}\rangle\\
&=d\mathscr{L}\times( d\Phi\times N)\\
&=d\mathscr{L}\times (*d\Phi).
\end{align}
Substituting these relations in the ones above then gives
\eqref{eq:Willmore-inh-system-2} and \eqref{eq:Willmore-inh-system-3}.
\end{bfproof}

\begin{proposition}
\label{prop:qualit-Lp-reg-Willmore}
Let $\Phi\in \mathscr{E}(B_1,\R^3)$
be conformal with conformal factor $\e^{\lambda}$
and so that $\delta\mathcal{W}\in L^p(B_1)$ for some $p>1$.
Then $H\in L^r_{\loc}(B_1)$ for every $r<\infty$.
\end{proposition}

\begin{bfproof}
We may certainly assume $1<p<2$.
\smallskip

\textit{Step 1: there exists
$\mathscr{L}$ and $L$ realizing the Hodge decomposition
\eqref{eq:Hodge-w} with $\mathscr{L}\in W^{2,p}(B_1)$
and $L\in L^{(2,\infty)}_{\loc}(B_1)$.} 
Indeed, we let $\mathscr{L}$ solve
\begin{align}
\left\{
\begin{aligned}
\Delta \mathscr{L} &= \e^{2\lambda} \delta\mathcal{W}
&&\text{in } B_1,\\
u&=0 &&\text{on }\partial B_1.
\end{aligned}
\right.
\end{align}
Elliptic regularity theory gives then that
$\mathscr{L}\in W^{2,p}$.
By construction we have
\begin{align}
d^*(w-d\mathscr{L})=0
\quad\text{in } B_1,
\end{align}
thus $L$ exists as a distribution in $B_1$ thanks to Poincar\'e's lemma and
it is determined up to an additive constant.
Note now that
\begin{align}
\Delta L
&= d^*(*d\mathscr{L}-*w)\\
&= d^*\big(*d\mathscr{L} 
- *(\nabla H 
+\langle A^\circ, H\rangle^{\sharp_g} 
+\langle A, H\rangle^{\sharp_g})\big)\\
&=-d^*\big(
*(\langle A^\circ, H\rangle^{\sharp_g} 
+\langle A, H\rangle^{\sharp_g})\big),
\end{align}
so if we let
\begin{align}
L_0(x)=-\int_{B_1}
\big\langle dK(x-y), *(\langle A^\circ, H\rangle^{\sharp_g} 
+\langle A, H\rangle^{\sharp_g})\big\rangle
\, dy,
\end{align}
where $K(x)=-\frac{1}{2\pi}\log|x|$ is the fundamental solution
of the Laplace operator,
given that
\begin{align}
|d K(x-y)|\leq \frac{C}{|x-y|}
\quad\Longrightarrow\quad
\sup_{x\in B_1}\|d K(x-\cdot)\|_{L^{2,\infty}(B_1)}
\leq C,
\end{align}
we see that $L_0\in L^{(2,\infty)}$,
and hence, since $L-L_0$ is harmonic, that
$L\in L^{(2,\infty)}_{\loc}.$
\smallskip

\textit{Step 2:
There exist $\mathscr{R}$, $R$, $\mathscr{S}$ and $S$ 
realizing the Hodge decompositions
\eqref{eq:Hodge-R}, \eqref{eq:Hodge-S} with
$\mathscr{R},\mathscr{S}\in W^{2,p^*}(B_1)$
and $R,S\in W^{1,(2,\infty)}_{\loc}(B_1)$.}
Indeed, define $\mathscr{R}$ and $\mathscr{S}$ by
\begin{align}
&\left\{
\begin{aligned}
\Delta \mathscr{R} &= -d\Phi\times d\mathscr{L}
&&\text{in } B_1,\\
\mathscr{R}&=0 &&\text{on }\partial B_1,
\end{aligned}
\right.
&
&\left\{
\begin{aligned}
\Delta \mathscr{S} &= -\langle d\Phi, d\mathscr{L}\rangle
&&\text{in } B_1,\\
\mathscr{S}&=0 &&\text{on }\partial B_1.
\end{aligned}
\right.
\end{align}
Since by Sobolev embedding
$d\mathscr{L}\in W^{1,p}\hookrightarrow L^{p^*}$,
so elliptic regularity gives
$\mathscr{R},\mathscr{S}\in W^{2,p^*}$.
By construction we then have
\begin{align}
d^*\big(-d\Phi\times H - (*d\Phi)\times L- d\mathscr{R}\big)&=0,\\
d^*\big(-\langle *d\Phi,L\rangle-d\mathscr{S} \big) &=0,
\end{align}
and thus $R$ and $S$ exist as distributions 
by Poincar\'e lemma and are determined up to 
additive constants,
and, since $L$ is in $L^{(2,\infty)}_{\loc}$,
so are $dR$ and $dS$.
\smallskip

\textit{Step 3: conclusion.}
From relations
\eqref{eq:Willmore-inh-system-2}, \eqref{eq:Willmore-inh-system-3},
we see that $R$ and $S$ satisfy a system Jacobians plus some extra terms,
namely
\begin{align}
\left\{
\begin{aligned}
\Delta R
&= dN\times(*dR) - \langle dN,*dS\rangle
+ f_R,\\
\Delta S
&= \langle dN,*dR\rangle
+ f_S,
\end{aligned}
\right.
\end{align}
where, since
$d\mathscr{L}\in L^{p^*}$ and
$d\mathscr{R}, d\mathscr{S}\in W^{1,p^*}\hookrightarrow L^\infty$,
we have
\begin{align}
f_R&=dN\times d\mathscr{R}
-\langle dN,d\mathscr{S}\rangle
+(*d\Phi)\times d\mathscr{L} &\in L^2,\\
f_S&=\langle dN,d\mathscr{R}\rangle
+\langle*d\Phi,d\mathscr{L}\rangle
&\in L^2.
\end{align}
Thanks to Theorem \ref{thm:Wente-variant},
we get that $R,S\in W^{2,q}_{\loc}$ for every $q<2$
and hence that 
\begin{align}
dR, dS\in W^{1,r}_{\loc} \quad\text{for every } r<\infty.
\end{align}
Inserting thin information in \eqref{eq:Willmore-inh-system-1}
gives that $\Delta\Phi$, and so $H$, is in 
$L^r_{\loc}$ for every $r<\infty$. 
\end{bfproof}

\subsection{Quantitative Estimates}
In the computations that follow, we
shall make use of various Gagliardo-Nirenberg inequalities,
namely, of multiplicative Sobolev inequalities such as
\begin{align*}
\|u\|_{L^4(\Omega)}\leq C
\|u\|_{L^2(\Omega)}^{1/2}
\|u\|_{W^{1,2}(\Omega)}^{1/2},
\end{align*}
for $u\in W^{1,2}(\Omega)$, $\Omega\subset\R^2$ bounded,
regular domain. 
We refer for instance to \cite{MR109940}.

\begin{proposition}
\label{prop:est-bilap-Phi-Willmore}
Let
$\Omega\subset\R^2$ is a bounded, regular domain
with $0\in\Omega$ and
$\Phi:\Omega\to \R^n$ is a conformal 
immersion of class $W^{4,2}$ with conformal factor $\e^{\lambda}$.
Let $E>0$ and $C_\infty>0$ be constants so that
\begin{align}
\|\e^{-\lambda}\nabla^2\Phi\|_{L^2(\Omega)}\leq E
\quad\text{and}\quad
\|\lambda - \Lambda\|_{L^\infty(\Omega)}
\leq C_{\infty},
\end{align}
where $\Lambda = \lambda(0)$.
Then, the following estimate holds:
\begin{align}
\label{eq:est-bilap-Phi-Willmore}
\|\Delta^2\Phi\|_{L^2(\Omega)}
\leq 
2\|\e^{4\lambda}\delta\mathcal{W}\|_{L^2(\Omega)}
+C \|\e^{-\lambda}\nabla^2\Phi\|_{L^2(\Omega)}
\|\nabla^2\Phi\|_{W^{2,2}(\Omega)},
\end{align}
for a constant
$C=C(\Omega, E, C_\infty)>0$.
\end{proposition}

First we point out a few elementary facts.
\begin{lemma}
For an immersion $\Phi:B_1\to\R^n$ there holds
\begin{align}
\label{eq:bilaplace-decomp}
\frac{1}{2}\Delta^2_g \Phi 
=\Delta_g H 
=\Delta^\perp_g H - \langle A,\langle H,A\rangle\rangle_g
-2\langle \nabla H, A\rangle^{\sharp_g}
-\langle \nabla H,H\rangle^{\sharp_g},
\end{align}
and, if $\Phi$ is conformal with conformal factor $\e^{\lambda}$,
there holds
\begin{align}
\label{eq:bilaplace-conformal}
\Delta^2 \Phi
&=\e^{4\lambda}\big(
\Delta_{g}^2\Phi + 2\langle d\lambda,\nabla H \rangle_g
+(2|d\lambda|^2_g +\Delta_{g}\lambda)H
\big).
\end{align}
\end{lemma}
The proof is elementary.

\begin{lemma}
The following pointwise estimates hold
for an absolute constant $C>0$ and $k=0,1,2$:
\begin{align}
\label{eq:est-hess-A}
\frac{1}{C}|\e^{-\lambda}\nabla^2\Phi|
&\leq |\e^{-\lambda}A| + |d\lambda|
\leq C\,\e^{-\lambda}|\nabla^2\Phi|,\\
\label{eq:est-H-Phi}
|\nabla^k(\e^{2\lambda}H)|
&\leq C|\nabla^{k+2}\Phi|,\\
\label{eq:est-H-A}
|\nabla^k(\e^{2\lambda}H)|
&\leq C|\nabla^{k}A|,
\end{align}
and the following estimate holds:
\begin{align}
\label{eq:est-A-hess}
\|A\|_{W^{k,2}(\Omega)}&\leq C\|\nabla^2\Phi\|_{W^{k,2}(\Omega)},
\end{align}
for $C=C(\Omega, E, C_\infty)>0$
and $k=0,1,2$.
\end{lemma}

\begin{bfproof}
Since
\begin{align}
A_{\mu\nu} 
= \partial^2_{\mu\nu}\Phi - \Gamma_{\mu\nu}^\sigma\,\partial_\sigma\Phi
\quad\text{and}\quad
e^{2\lambda}=\frac{1}{2}|d\Phi|^2,
\end{align}
and
\begin{align}
\Gamma_{11}^1&=\partial_1\lambda,&
\Gamma^1_{12}&=\Gamma^1_{21}=\partial_2\lambda,&
\Gamma^1_{22}&=-\partial_1\lambda,\\
\Gamma_{11}^2&=-\partial_2\lambda,&
\Gamma^2_{12}&=\Gamma^1_{21}=\partial_1\lambda,&
\Gamma^2_{22}&=\partial_2\lambda,
\end{align}
estimates \eqref{eq:est-hess-A} and consequently 
\eqref{eq:est-A-hess} for $k=0$ follow.
Next, since
\begin{align}
\e^{2\lambda}H
&=\frac{1}{2}\Delta\Phi
=\frac{1}{2}(\Delta\Phi)^\perp
=\frac{1}{2} (A_{11}+A_{22}),
\end{align}
estimates \eqref{eq:est-H-Phi} and \eqref{eq:est-H-A} follow.
Now differentiating of the above identities:
\begin{align}
\partial_\xi A_{\mu\nu}
&=\partial^3_{\xi\mu\nu}\Phi 
- \partial_\xi\Gamma_{\mu\nu}^\sigma\,\partial_\sigma\Phi
- \Gamma_{\mu\nu}^\sigma\,\partial^2_{\xi\sigma}\Phi,\\
\partial^2_{\zeta\xi} A_{\mu\nu}
&=\partial^4_{\zeta\xi\mu\nu}\Phi 
- \partial^2_{\zeta\xi}\Gamma_{\mu\nu}^\sigma\,\partial_\sigma\Phi
- \partial_{\xi}\Gamma_{\mu\nu}^\sigma\,\partial^2_{\zeta\sigma}\Phi
- \partial_\zeta\Gamma_{\mu\nu}^\sigma\,\partial^2_{\xi\sigma}\Phi
- \Gamma_{\mu\nu}^\sigma\,\partial^3_{\zeta\xi\sigma}\Phi,
\end{align}
and
\begin{align}
2d\lambda\e^{2\lambda} &=\langle \nabla^2\Phi,d\Phi\rangle,\\
(2\nabla^2\lambda +4d\lambda\otimes d\lambda)\e^{2\lambda}
&=\langle\nabla^3\Phi,d\Phi\rangle
+\nabla^2\Phi\langle\dot{\otimes}\rangle\nabla^2\Phi,\\
(2\nabla^3\lambda +8\nabla^2\lambda\otimes d\lambda
+8d\lambda\otimes d\lambda \otimes d\lambda)\e^{2\lambda}
&=\langle\nabla^4\Phi,d\Phi\rangle
+2\nabla^3\Phi\langle \dot{\otimes}\rangle \nabla^2\Phi,
\end{align}
(where $\langle\dot{\otimes}\rangle$ means
scalar product in the vector part, inner product in one
of the covariant entries and tensor product in the remaining ones)
yields the estimates
\begin{align}
|\Gamma_{\mu\nu}^\sigma|&\leq C |d\lambda|,&
|\partial_\xi\Gamma_{\mu\nu}^\sigma|&\leq C |\nabla^2\lambda|,&
|\partial^2_{\zeta\xi}\Gamma_{\mu\nu}^\sigma|&\leq C |\nabla^3\lambda|,
\end{align}
and
\begin{align}
|d\lambda|&\leq C\e^{-\lambda}|\nabla^2\Phi|,\\
|\nabla^2\lambda|
&\leq C\big(
|d\lambda|^2 + \e^{-\lambda}|\nabla^3\Phi| + \e^{-2\lambda}|\nabla^2\Phi|^2
\big)\\
&\leq C\big(
\e^{-2\lambda}|\nabla^2\Phi|^2 + \e^{-\lambda}|\nabla^3\Phi|
\big),\\
|\nabla^3\lambda|
&\leq C\big(
|\nabla^2\lambda||d\lambda| + |d\lambda|^3
+\e^{-\lambda}|\nabla^4\Phi|
+\e^{-2\lambda}|\nabla^3\Phi||\nabla^2\Phi|
\big)\\
&\leq C\big(
\e^{-2\lambda}|\nabla^3\Phi||\nabla^2\Phi|
+\e^{-3\lambda}|\nabla^2\Phi|^3
+\e^{-\lambda}|\nabla^4\Phi|
\big),
\end{align}
thanks to which we estimate in turn
\begin{align}
|\nabla A|
&\leq C\big(
|\nabla^3\Phi| + \e^{\lambda}|\nabla^2\lambda| + |d\lambda||\nabla^2\Phi|
\big)\\
&\leq C\big(
|\nabla^3\Phi| + \e^{-\lambda}|\nabla^2\Phi|^2
\big),\\
|\nabla^2A|
&\leq C\big(
|\nabla^4\Phi| + \e^{\lambda}|\nabla^3\lambda|
+|\nabla\lambda||\nabla^2\Phi| +|d\lambda||\nabla^3\Phi|
\big)\\
&\leq C\big(|\nabla^4\Phi| + \e^{-\lambda}|\nabla^3\Phi||\nabla^2\Phi|
+\e^{-2\lambda}|\nabla^2\Phi|^3
\big).
\end{align}
Thus, with the help of Gagliardo-Nirenberg we can estimate
\begin{align}
\|\nabla A\|_{L^2}
&\leq C\big(
\|\nabla^3\Phi\|_{L^2}
+\|\e^{-\lambda}|\nabla^2\Phi|^2\|_{L^2}
\big)\\
&\leq C\big(
\|\nabla^3\Phi\|_{L^2}
+\e^{-\Lambda}\|\nabla^2\Phi\|_{L^4}^2
\big)\\
&\leq C\big(
\|\nabla^3\Phi\|_{L^2}
+\e^{-\Lambda}\|\nabla^2\Phi\|_{L^2}\|\nabla^2\Phi\|_{W^{1,2}}
\big)\\
&\leq C\big(
\|\nabla^3\Phi\|_{L^2}
+\|\e^{-\lambda}\nabla^2\Phi\|_{L^2}\|\nabla^2\Phi\|_{W^{1,2}}
\big),
\end{align}
yielding \eqref{eq:est-A-hess} for $k=1$.
Similarly, since
\begin{align}
\|\nabla^2A\|_{L^2}
&\leq C\big(
\|\nabla^4\Phi\|_{L^2}
+ \|\e^{-\lambda}|\nabla^3\Phi||\nabla^2\Phi|\|_{L^2}
+\|\e^{-2\lambda}|\nabla^2\Phi|^3\|_{L^2}
\big)\\
&\leq C\big(
\|\nabla^4\Phi\|_{L^2}
+ \e^{-\Lambda}\||\nabla^3\Phi||\nabla^2\Phi|\|_{L^2}
+\e^{-2\Lambda}\|\nabla^2\Phi\|_{L^6}^3
\big),
\end{align}
with H\"older and Gagliardo-Nirenberg we estimate
\begin{align}
\||\nabla^3\Phi||\nabla^2\Phi|\|_{L^2}
&\leq\|\nabla^3\Phi\|_{L^4}\|\nabla^2\Phi\|_{L^4}\\
&\leq C\|\nabla^2\Phi\|_{L^2}^{\frac{1}{4}}\|\nabla^2\Phi\|_{W^{2,2}}^{\frac{3}{4}}
\|\nabla^2\Phi\|_{L^2}^{\frac{3}{4}}\|\nabla^2\Phi\|_{W^{2,2}}^{\frac{1}{4}}\\
&\leq C \|\nabla^2\Phi\|_{L^2}\|\nabla^2\Phi\|_{W^{2,2}},
\end{align}
and
\begin{align}
\|\nabla^2\Phi\|_{L^6}^3
\leq C\|\nabla^2\Phi\|_{L^2}^2\|\nabla^2\Phi\|_{W^{2,2}},
\end{align}
so to obtain
\begin{align}
\|\nabla^2A\|_{L^2}
&\leq C\big(
\|\nabla^4\Phi\|_{L^2}
+\|\e^{-\lambda}\nabla^2\Phi\|_{L^2}\|\nabla^2\Phi\|_{W^{2,2}}
+\|\e^{-\lambda}\nabla^2\Phi\|_{L^2}^2\|\nabla^2\Phi\|_{W^{2,2}}
\big),
\end{align}
yielding \eqref{eq:est-A-hess} also for $k=2$.
\end{bfproof}

\begin{bfproof}[Proof of Proposition \ref{prop:est-bilap-Phi-Willmore}]
It can be deduced from the following three lemmas
and \eqref{eq:est-A-hess}.
\begin{lemma}
\label{lemma:est-normlap-H-A}
There holds
\begin{align}
\label{eq:est-normlap-H-A}
\|\e^{4\lambda}\Delta^\perp_g H\|_{L^2(\Omega)}
\leq C
\|\e^{-\lambda}A\|_{L^2(\Omega)}\|A\|_{W^{2,2}(\Omega)}
+\|\e^{4\lambda}\delta\mathcal{W}\|_{L^2(\Omega)},
\end{align}
for a constant $C=C(\Omega, E, C_\infty)>0$.
\end{lemma}

\begin{lemma}
\label{lemma:est-lap-H-normlap-H}
There holds
\begin{align}
\label{eq:est-lap-H-normlap-H}
\|\e^{4\lambda}\Delta_g H\|_{L^2(\Omega)}\leq
\|\e^{4\lambda}\Delta^\perp_g H\|_{L^2(\Omega)}
+C\|\e^{-\lambda}A\|_{L^2(\Omega)}\|A\|_{W^{2,2}(\Omega)},
\end{align}
for a constant $C=C(\Omega, E, C_\infty)>0$.
\end{lemma}

\begin{lemma}
\label{lemma:est-bilap-phi-lap-H}
There holds
\begin{align}
\label{eq:est-bilap-phi-lap-H}
\|\Delta^2\Phi\|_{L^2(\Omega)}
\leq 2\|\e^{4\lambda}\Delta_g H\|_{L^2(\Omega)} 
+C \|\e^{-\lambda}\nabla^2\Phi\|_{L^2(\Omega)}
\|\nabla^2\Phi\|_{W^{2,2}(\Omega)},
\end{align}
for a constant
$C=C(\Omega,E, C_\infty)>0$.
\end{lemma}

\begin{bfproof}
[Proof of Lemma \ref{lemma:est-normlap-H-A}]
With the classical form of the Willmore operator \eqref{eq:Willmore-classical}
and the pointwise estimate \eqref{eq:est-H-A},
by means of H\"older's and Gagliardo-Nirenberg inequalities we estimate
\begin{align}
\|\langle A^\circ,\langle H, A^\circ\rangle\rangle_g\|_{L^2}
&\leq C
\|\e^{-4\lambda}\langle A^\circ,\langle H, A^\circ\rangle\rangle\|_{L^2}\\
&\leq C\|\e^{-4\lambda}|A^\circ|^2 H\|_{L^2}\\
&\leq C\|\e^{-6\lambda}|A|^3\|_{L^2}\\
&\leq C\e^{-6\Lambda}\|A\|_{L^6}^3\\
&\leq C\e^{-6\Lambda}\|A\|_{L^2}^2\|A\|_{W^{2,2}}\\
&\leq  C\e^{-4\Lambda}\|\e^{-\lambda}A\|_{L^2}^2\|A\|_{W^{2,2}},
\end{align}
which yields \eqref{eq:est-normlap-H-A}.
\end{bfproof}

\begin{bfproof}
[Proof of Lemma \ref{lemma:est-lap-H-normlap-H}]
From formula \eqref{eq:bilaplace-decomp},
we have
\begin{align}
\|\Delta_g H\|_{L^2}
\leq \|\Delta^\perp_g H\|_{L^2}
+C\big(
\|\langle A,\langle H, A\rangle\rangle_g\|_{L^2}
+\|\langle \nabla H, A\rangle^{\sharp_g}_g\|_{L^2}
+\|\langle \nabla H, H\rangle^{\sharp_g}\|_{L^2}
\big).
\end{align}
Similarly as in the proof of Lemma \eqref{lemma:est-normlap-H-A}, we
have
\begin{align}
\|\langle A,\langle H, A\rangle\rangle_g\|_{L^2}
\leq C\e^{-4\Lambda}\|\e^{-\lambda}A\|_{L^2}^2\|A\|_{W^{2,2}}.
\end{align}
Next, since
\begin{align}
\nabla H
&= \e^{-2\lambda}\nabla(\e^{2\lambda} H)
-2H\otimes d\lambda,
\end{align}
with \eqref{eq:est-H-A} we may pointwise estimate
\begin{align}
|\nabla H|&\leq C\e^{-2\lambda}
\big(|\nabla A| + |A||d\lambda|\big),
\end{align}
thus allowing to deduce
\begin{align}
\|\langle \nabla H, H\rangle^{\sharp_g}\|_{L^2}
&\leq C\|\e^{-\lambda}|\nabla H| |H|\|_{L^2}\\
&\leq C\|\e^{-3\lambda}(|\nabla A| + |A||d\lambda|)|H|\|_{L^2}\\
&\leq C\|\e^{-5\lambda}(|\nabla A| + |A||d\lambda|)|A|\|_{L^2};
\end{align}
now with H\"older and Gagliardo-Nirenberg we see that,
on the one hand,
\begin{align}
\|\e^{-5\lambda}|\nabla A||A|\|_{L^2}
&\leq C \e^{-5\Lambda}\||\nabla A||A|\|_{L^2}\\
&\leq C \e^{-5\Lambda}\|\nabla A\|_{L^4}\|A\|_{L^4}\\
&\leq C \e^{-5\Lambda}\| A\|_{L^2}\|A\|_{W^{2,2}}\\
&\leq C \e^{-4\Lambda}\|\e^{-\lambda} A\|_{L^2}\|A\|_{W^{2,2}},
\end{align}
and on the other hand
\begin{align}
\|\e^{-5\lambda}|A|^2|d\lambda|\|_{L^2}
&\leq C\e^{-5\Lambda}\|d\lambda\|_{L^2}
\|A\|_{L^\infty}^2\\
&\leq C\e^{-5\Lambda}\|d\lambda\|_{L^2}
\|A\|_{L^2}\|A\|_{W^{2,2}}\\
&\leq C\e^{-4\Lambda}\|d\lambda\|_{L^2}
\|\e^{-\lambda}A\|_{L^2}\|A\|_{W^{2,2}},
\end{align}
so we deduce
\begin{align}
\|\langle \nabla H, H\rangle^{\sharp_g}\|_{L^2}
\leq C\e^{-4\Lambda}
\|\e^{-\lambda}A\|_{L^2}\|A\|_{W^{2,2}}.
\end{align}
Similarly, we see that
\begin{align}
\|\langle \nabla H, A\rangle^{\sharp_g}_g\|_{L^2}
&\leq C\|\e^{-3\lambda}\langle \nabla H, A\rangle^{\sharp}\|_{L^2}\\
&\leq C\|\e^{-3\lambda}|\nabla H| |A|\|_{L^2}\\
&\leq C\|\e^{-5\lambda}(|\nabla A| + |A||d\lambda|) |A|\|_{L^2},
\end{align}
and so similarly as before we deduce
\begin{align}
\|\langle \nabla H, A\rangle^{\sharp_g}_g\|_{L^2}
\leq C\e^{-4\Lambda}\|\e^{-\lambda}A\|_{L^2}\|A\|_{W^{2,2}},
\end{align}
yielding \eqref{eq:est-lap-H-normlap-H}.
\end{bfproof}

\begin{bfproof}
[Proof of Lemma \ref{lemma:est-bilap-phi-lap-H}]
From formula \eqref{eq:bilaplace-conformal}, it follows that
\begin{align}
\|\Delta\Phi\|_{L^2}
\leq
2\|\e^{4\lambda}\Delta_g H\|_{L^2}
+ C\big(
\|\e^{4\lambda}\langle d\lambda,\nabla H\rangle_g\|_{L^2}
+\|\e^{4\lambda}|d\lambda|^2_g H\|_{L^2}
+ \|\e^{4\lambda}\Delta_g \lambda\,H\|_{L^2}
\big).
\end{align}
Now with \eqref{eq:est-hess-A}, \eqref{eq:est-H-Phi}, and (again) the identity
\begin{align}
\nabla H
&= \e^{-2\lambda}\nabla(\e^{2\lambda} H)
-2H\otimes d\lambda,
\end{align}
we estimate with H\"older and Galgliardo-Nirenberg
\begin{align}
\|\e^{4\lambda}\langle d\lambda,\nabla H\rangle_g\|_{L^2}
&\leq C\|\e^{2\lambda} \langle d\lambda,\nabla H\rangle\|_{L^2}\\
&\leq C\|\e^{2\lambda} |d\lambda||\nabla H|\|_{L^2}\\
&\leq C\big\|\e^{2\lambda}\big(\e^{-\lambda}|\nabla^2\Phi|\big)
\big(\e^{-2\lambda}|\nabla^3\Phi| + \e^{-3\lambda}|\nabla^2\Phi|^2\big)\big\|_{L^2}\\
&\leq C\big(\e^{-\Lambda}\||\nabla^3\Phi|\nabla^2\Phi|\|_{L^2}
+\e^{-2\Lambda}\|\nabla^2\Phi\|_{L^6}^3\big)\\
&\leq C\big(\e^{-\Lambda}\|\nabla^3\Phi\|_{L^4}\|\nabla^2\Phi\|_{L^4}
+\e^{-2\Lambda}\|\nabla^2\Phi\|_{L^6}^3\big)\\
&\leq C\big(\e^{-\Lambda}\|\nabla^2\Phi\|_{L^2}\|\nabla^2\Phi\|_{W^{2,2}}
+\e^{-2\Lambda}\|\nabla^2\Phi\|_{L^2}^2\|\nabla^2\Phi\|_{W^{2,2}}\big)\\
&\leq C\|\e^{-\lambda}\nabla^2\Phi\|_{L^2}\|\nabla^2\Phi\|_{W^{2,2}}.
\end{align}
Similarly, again with \eqref{eq:est-H-Phi}
and Gagliardo-Nirenberg we estimate
\begin{align}
\|\e^{4\lambda}|d\lambda|^2_g H\|_{L^2}
&\leq C\|\e^{2\lambda}|d\lambda|^2 H\|_{L^2}\\
&\leq C\|\e^{2\lambda}(\e^{-2\lambda}|\nabla^2\Phi|^2)
(\e^{-2\lambda}|\nabla^2\Phi|)\|_{L^2}\\
&\leq C\|\e^{-2\lambda}|\nabla^2\Phi|^3\|_{L^2}\\
&\leq C\|\e^{-\lambda}\nabla^2\Phi\|_{L^2}
\|\nabla^2\Phi\|_{W^{2,2}}.
\end{align}
Finally, with Liouville's equation
\begin{align}
-\Delta\lambda =\e^{2\lambda}K, 
\end{align}
the pointwise estimate
\begin{align}
|K|\leq |A|^2_g\leq \e^{-4\lambda}|A|^2,
\end{align}
we can estimate with Gagliardo-Nirenberg
and \eqref{eq:est-A-hess}:
\begin{align}
\|\e^{4\lambda}\Delta_g\lambda H\|_{L^2}
&\leq C\||A|^2 H\|_{L^2}\\
&\leq C\|\e^{-2\lambda}|A|^3\|_{L^2}\\
&\leq C\e^{-2\Lambda}\|A\|_{L^2}^2\|A\|_{W^{2,2}}\\
&\leq C\|\e^{-\lambda}A\|_{L^2}^2\|A\|_{W^{2,2}}\\
&\leq C\|\e^{-\lambda}A\|_{L^2}^2\|\nabla^2\Phi\|_{W^{2,2}}.
\end{align}
All these estimates together yield \eqref{eq:est-bilap-phi-lap-H}.
\end{bfproof}
The combination of lemmas \ref{lemma:est-normlap-H-A},
\ref{lemma:est-lap-H-normlap-H} and 
\ref{lemma:est-bilap-phi-lap-H}
immediately gives \eqref{eq:est-bilap-Phi-Willmore},
and concludes the proof of Proposition \ref{prop:est-bilap-Phi-Willmore}.
\end{bfproof}

\begin{bfproof}
[Proof of Theorem \ref{thm:quant-L2-reg-Willmore}]
The qualitative statement, namely that $\Phi\in W^{4,2}_{\loc}$,
is an immediate consequence of Propositions \ref{prop:bootstrap-cond-Willmore}
and \ref{prop:qualit-Lp-reg-Willmore}.
We can therefore apply Proposition \ref{prop:est-bilap-Phi-Willmore},
so that
\eqref{eq:est-bilap-Phi-Willmore} jointly with
elliptic estimates for the bilaplacian give
\begin{align}
\label{eq:quant-L2-est-1st-step}
\|d\Phi\|_{W^{3,2}(B_{1/2})}
\leq C\big(
\|\e^{4\lambda}\delta\mathcal{W}\|_{L^2(B_1)}
+\|\e^{-\lambda}\nabla^2\Phi\|_{L^2(B_1)}\|\nabla^2\Phi\|_{W^{2,2}(B_1)}
+\|d\Phi\|_{L^2(B_1)}
\big)
\end{align}
for $C=C(E, C_{\infty})>0$.
Now we consider rescalings.
For $0<r<1$ we let
\begin{align}
\widetilde{\Phi}(x) = \Phi(rx),
\quad x\in B_1,
\end{align}
and we denote with a tilde all the quantities 
pertaining to $\widetilde{\Phi}$.
From
for $k\in \N$ it follows in particular that
for $\Omega\subseteq B_1$ we have
\begin{align}
\nabla^k\widetilde{\Phi}(x) &= r^k\nabla^k\Phi(rx),\\
\|\nabla^k\widetilde{\Phi}\|_{L^2(\Omega)}
&=r^{k-1}\|\nabla^k\Phi\|_{L^2(r\Omega)},\\
\e^{\widetilde{\lambda}(x)} 
&=r\e^{\lambda(rx)},\\
\widetilde{\lambda}(x) - \widetilde{\lambda}(0)
&=\lambda(rx) - \lambda(0),\\
\widetilde{\delta\mathcal{W}}(x)
&=\delta\mathcal{W}(rx),
\end{align}
(the last equality follows either by direct inspection
or at once recalling that $\delta\mathcal{W}$ is a vector field)
and from these relations, we deduce in particular that
\begin{align}
\|\e^{4\widetilde{\lambda}}\widetilde{\delta\mathcal{W}}\|_{L^2(B_1)}
&=r^3\|\e^{4\lambda}\delta\mathcal{W}\|_{L^2(B_r)},\\
\|\e^{-2\widetilde{\lambda}}\nabla^2\widetilde{\Phi}\|_{L^2(B_1)}
&=\|\e^{-2\lambda}\nabla^2\Phi\|_{L^2(B_r)},\\
\|\widetilde{\lambda} - \widetilde{\Lambda}\|_{L^\infty(B_1)}
&= \|\lambda - \Lambda\|_{L^\infty(B_r)},
\end{align}
and  the last two relations in particular give that
$\widetilde{E}\leq E$ and $\widetilde{C_{\infty}}\leq C_{\infty}$.
Consequently,
applying \eqref{eq:quant-L2-est-1st-step} to
$\widetilde{\Phi}$ gives
{\small
\begin{align}
\label{eq:est-invariant-sob}
\|d\Phi\|'_{W^{3,2}(B_{r/2})}
&\leq C\big(
r^3\|\e^{4\lambda}\delta\mathcal{W}\|_{L^2(B_r)}
+\|\e^{-\lambda}\nabla^2\Phi\|_{L^2(B_r)}\|\nabla^2\Phi\|_{W^{2,2}(B_r)}'
+\|d\Phi\|_{L^2(B_r)}
\big),
\end{align}
}
where for $k=1,2$ we denoted
$
\|\Phi\|_{W^{k,2}(B_{\rho})}'
=\big(\sum_{h=1}^k\rho^{k-1}\|\Phi\|_{W^{h,2}(B_{\rho})}^2\big)^{1/2}
$
the scale-invariant version of the Sobolev norms
where, as we said abefore $C= C(E,C_\infty)>0$.
\medskip 

Recall now that from \eqref{eq:est-hess-A}
it is
\begin{align}
\|\e^{-\lambda}\nabla^2\Phi\|_{L^2(B_{1/2})}
\leq C\big(
\|\e^{-\lambda}A\|_{L^2(B_{1/2})}
+\|d\lambda\|_{L^2(B_{1/2})}
\big),
\end{align}
so by Theorem \ref{thm:eps-control-conformal-fact}
we let $\varepsilon_0$ be sufficiently small so to have
\begin{align}
\|d\lambda\|_{L^2(B_{1/2})}
+\|\lambda - \lambda(0)\|_{L^\infty(B_{1/2})}
\leq C\|\e^{-\lambda}A\|_{L^2(B_{1})}
\leq C\varepsilon_0,
\end{align}
for  $C=C(C_{(2,\infty)})>0$,
whence \eqref{eq:est-invariant-sob} can be improved to
\begin{align}
\|d\Phi\|'_{W^{3,2}(B_{r/2})}
&\leq C\big(
r^3\|\e^{4\lambda}\delta\mathcal{W}\|_{L^2(B_r)}
+\varepsilon_0\|\nabla^2\Phi\|_{W^{2,2}(B_r)}'
+\|d\Phi\|_{L^2(B_r)}
\big),
\end{align}
for $C=C(C_{(2,\infty)})>0.$
\medskip

Choose finally $\varepsilon_0>0$ be sufficiently 
small so that $C\varepsilon_0\leq\frac{1}{2}$
to obtain
\begin{align}
\|d\Phi\|'_{W^{3,2}(B_{r/2})}
&\leq
\frac{1}{2}\|d\Phi\|'_{W^{3,2}(B_{r})}
+C\big(
r^3\|\e^{4\lambda}\delta\mathcal{W}\|_{L^2(B_1)}
+\|d\Phi\|_{L^2(B_1)}
\big),
\end{align}
for every $0<r\leq 1/2$.
A classical iteration/interpolation argument applied
to $\phi(r) = \|d\Phi\|'_{W^{3,2}(B_{r})}$ and a covering
argument yield to
\eqref{eq:quant-L2-est-Willmore}.
\end{bfproof}

\section{Conformal Willmore Flows}
\label{sec:conf-flows}
\begin{bfproof}[Proof of Lemma \ref{lemma:conformal-flow-cond}]
A metric $g$ is conformal if and only if
its Hopf differential (computed with respect to the background complex
structure of $S^2$) vanishes identically. In our case it is
\begin{align}
\label{eq:complex-bil-part-metric}
\operatorname{Hopf}(g)
=g_{zz}dz\otimes dz
=\langle\partial_z\Phi,\partial_z\Phi\rangle
dz\otimes dz.
\end{align}
Since $\delta\mathcal{W}$ is a normal vector field,
we see that
\begin{align}
\frac{1}{2}\frac{\partial}{\partial t}\langle\partial_z\Phi,\partial_z\Phi\rangle
&=
\langle\partial_z\partial_t\Phi,\partial_z\Phi\rangle
=
\langle
-\partial_z \delta\mathcal{W}+
\partial_z U,\partial_z\Phi
\rangle
=
\langle\delta\mathcal{W},
\partial^2_{zz}\Phi\rangle
+\langle\partial_z U,
\partial_z\Phi\rangle.
\end{align}
Since 
$U= U^{z}\partial_z\Phi + U^{\bar{z}}\partial_{\bar{z}}\Phi$
with $U^{z}=U^1+iU^2$ and $U^{\bar{z}}=U^1-iU^2$
we have
\begin{align}
\left\langle\partial_z U,\partial_z\Phi\right\rangle
&=\big\langle
\partial_z U^z\partial_z\Phi 
+ U^{z}\partial^2_{zz}\Phi
+\partial_{z}U^{\bar{z}}\partial_{\bar{z}}\Phi
+U^{\bar{z}}\partial^2_{z\bar{z}}\Phi,
\partial_z\Phi
\big\rangle\\
&=\partial_z U^{z}g_{zz}
+U^{z}
\big\langle
\partial^2_{zz}\Phi,\partial_z\Phi
\big\rangle
+\partial_zU^{\bar{z}}g_{\bar{z}z}
+U^{\bar{z}}
\big\langle
\partial^2_{z\bar{z}}\Phi,\partial_z\Phi
\big\rangle\\
&=\partial_zU^{z}g_{zz} +
\partial_zU^{\bar{z}}g_{\bar{z}z}
+\frac{1}{2}
\Big(
U^z\partial_zg_{zz} + U^{\bar{z}}\partial_{\bar{z}}g_{zz}
\Big),
\end{align}
thus
we have
\begin{align}
\frac{1}{2}\frac{\partial}{\partial t}\langle\partial_z\Phi,\partial_z\Phi\rangle
&=
\partial_zU^{z}g_{zz} +
\partial_zU^{\bar{z}}g_{\bar{z}z}
+\frac{1}{2}
\Big(
U^z\partial_zg_{zz} + U^{\bar{z}}\partial_{\bar{z}}g_{zz}
\Big)
+\big\langle \delta\mathcal{W},\partial^2_{zz}\Phi
\big\rangle.
\end{align}
If the flow $\Phi$ is conformal,
then 
$\partial_t\operatorname{Hopf}(g)\equiv 0$
and
$g_{zz}, g_{\bar{z}\bar{z}}$ vanish identically. 
Since moreover $\delta\mathcal{W}$ is a normal vector we
may replace $\partial^2_{zz}\Phi$ with $A_{zz}$ and thus 
obtain, after conjugation,
\begin{align}
g_{\bar{z}z}\partial_z U^{\bar{z}} 
= -\langle\delta\mathcal{W},A_{zz}\rangle,
\end{align}
and so since $\frac{1}{2}\e^{2\lambda}=g_{\bar{z}{z}}$,
this yields \eqref{eq:conformal-flow-cond}.
\end{bfproof}

The proof of the following proposition follows directly by
direct inspection of the proof of the original theorems;
for the proof of (ii) one uses additionally that conformal
Willmore flows are $W^{4,2}$ for almost every time, as proved
in Proposition \ref{prop:est-conform-flow} below.

\begin{proposition}
\label{prop:nonsmooth-area-DLM}
\begin{enumerate}[label=(\roman*)]
\item If the surface $\mathcal{S}$ is immersed
through a $W^{4,2}$-map,
then the same conclusion of Theorem \ref{thm:DeLellis-Muller}
and all its consequences obtained in Section \ref{sec:DLM-consequences}
still hold.
\item For conformal Willmore flows in $\mathscr{W}^{\varepsilon,\delta}_{[0,T)}(S^2,\R^3)$
the area and barycenter bounds in Theorem \ref{thm:area-est}
still hold if $\varepsilon$ is chosen sufficiently small.
\end{enumerate}
\end{proposition}
A first regularity improvement for conformal Willmore flows
is the following.
\begin{proposition}
\label{prop:est-conform-flow}
For any $\varepsilon_0,\,\delta>0$,
any conformal Willmore flow 
$\Phi\in\mathscr{W}_{[0,T)}^{\varepsilon_0,\delta}(S^2,\R^3)$
is in $W^{4,2}(S^2)$ for a.e. $t\in (0,T)$.
Futhermore there exist $\varepsilon_0,\,\delta,\,C>0$
independent of $\Phi$ so that
$\Phi\in L^2((0,T),W^{4,2}(S^2))$ with
\begin{align}
\label{eq:est-conform-flow}
\|d\Phi\|_{L^2((0,T),W^{3,2}(S^2))}
\leq C(\sqrt{T}+
\|\e{^\lambda}\delta\mathcal{W}\|_{L^2((0,T)\times S^2)}),
\end{align}
for a constant $C>0$.
\end{proposition}

\begin{bfproof}
Since $\delta\mathcal{W}\in L^2((0,T)\times S^2)$,
Fubini's theorem implies $\delta\mathcal{W}(t,\cdot)\in L^2(S^2)$
for a.e. $t$. 
Consequently, by Propositions \ref{prop:bootstrap-cond-Willmore}
and \ref{prop:qualit-Lp-reg-Willmore} it follows that $\Phi(t,\cdot)\in  W^{4,2}(S^2)$.

From Liouville's equation
\begin{align}
\Delta_{S^2}\lambda =\e^{2\lambda}K - 1,
\end{align}
the pointwise inequality $|K|\leq |A|^2_g=\e^{-4\lambda}|A|^2$
and elliptic estimates
\begin{align}
\label{eq:est-conf-flow-1}
\|d\lambda\|_{L^\infty((0,T),L^{(2,\infty)}(S^2))}
\leq C.
\end{align}
Choose now $\varepsilon_1>0$ sufficiently small so that,
for a.e. $t\in(0,T)$ so that $\Phi(t,\cdot)$ is $W^{4,2}$,
we can fix a value $r>0$ which satisfies,
according to
Lemma \ref{lemma:DeLellis-Muller-local},
\begin{align}
\label{eq:est-conf-flow-2}
\int_{B_r(x_0)}
|A|^2\,d\sigma
&\leq C\Big(\varepsilon_1+
\e^{C}r^2\Big)\leq \varepsilon_0,
\end{align}
for every $x_0\in S^2$,
where $\varepsilon_0$ is as in Theorem \ref{thm:quant-L2-reg-Willmore}.
With the estimates \eqref{eq:est-conf-flow-1}
and \eqref{eq:est-conf-flow-2},
an application of Theorem \ref{thm:quant-L2-reg-Willmore}
to a covering $S^2$ with balls of radius $r/2$
gives, for a.e. $t\in(0,T)$,
\begin{align}
\|d\Phi\|_{W^{3,2}(S^2)}
\leq C\big(\|\e^{4\lambda}\delta\mathcal{W}\|_{L^2(S^2)}
+\|\e^{\lambda}\|_{L^2(S^2)}\big)
\leq C\big(\|\e^{\lambda}\delta\mathcal{W}\|_{L^2(S^2)}
+1\big),
\end{align}
where the last inequality follows recalling that by
assumption the bound \eqref{eq:def-class-DLM}
holds in the class $\mathscr{W}_{[0,T)}^{\varepsilon_0,\delta}(S^2,\R^3)$.
Now if we square and integrate in $t$ such inequality,
recalling that, since 
$\Phi\in\mathscr{W}_{[0,T)}^{\varepsilon_0,\delta}$,
\eqref{eq:weaker-energy-id} holds by assumption,
we obtain \eqref{eq:est-conform-flow}.
\medskip
\end{bfproof}

One can see that, in fact,  along
the proof of Proposition \ref{prop:est-conform-flow},
well-balanced condition (iv) in Definition
\ref{def:energy-class} was not needed,
an in fact the proof is completely independent
of the tangential component $U$.
Such condition is however essential to prove the following.

\begin{proposition}
\label{prop:est-U}
For every $p<2$ there exist $\delta>0$
with the following property.
Let $\Phi\in \mathscr{W}^{\varepsilon,\delta}_{[0,T)}(S^2,\R^3)$
be a weak Willmore flow.
Then its tangential component
is in $L^2((0,T),L^p(S^2))$,
with
\begin{align}
\label{eq:est-U}
\|U\|_{L^2((0,T),L^p(S^2))}
\leq C
\|\e^{\lambda}\delta\mathcal{W}\|_{L^2((0,T)\times S^2))},
\end{align}
for a constant $C=C(p)>0$.
\end{proposition}

\begin{bfproof}
We may certainly assume $p>1$.
In this proof we find it convenient to clearly distinguish
between the non-immersed (pulled-back) tangential component 
$U=U^\mu\partial_\mu$,
the immersed one $d\Phi(U)=U^\mu\partial_\mu\Phi$,
and the associated tangent vector field
on $S^2$, $dI(U)=U^\mu\partial_\mu I$,
where $I:S^2\to\R^3$ denotes is the standard immersion.
With this notation $\Phi$ satisfies weakly
\begin{align*}
\frac{\partial}{\partial t}\Phi = -\delta\mathcal{W} + d\Phi(U)
\quad\text{in }(0,T)\times S^2.
\end{align*}
In what follows, it will be implicitly understood
that all the  slice-wise operations are 
valid for a.e. fixed $t$.
Finally we shall make use of the notation and
concepts recalled in  Appendix \ref{sec:delbar}.
\medskip

From Lemma \ref{lemma:conformal-flow-cond},
we deduce that $U$ is given by
\begin{align}
U^{(1,0)}=
-\overline{\partial}{}^{-1}
(\langle\delta\mathcal{W},\overline{h}_0\rangle^{\sharp_g}) + \Omega
\end{align}
for some time-dependent holomorphic vector field $\Omega=\Omega(t,\cdot)\in\mathfrak{X}^{\omega}(S^2)$.
Classical elliptic estimates 
and H\"older's inequality permit to estimate
\begin{align}
\label{eq:est-U-ellipt}
\|\overline{\partial}{}^{-1}
(\langle\delta\mathcal{W},\overline{h}_0\rangle^{\sharp_g})
\|_{L^p(S^2)}
&\leq C_p\|\langle \delta\mathcal{W},
\overline{h}_0\rangle^\sharp_g\|_{L^1(S^2)}\\
&\leq C_p\||\delta\mathcal{W}|\e^{-2\lambda}|A^\circ|\|_{L^1(S^2)}\\
&\leq C_p\|\e^{-\lambda}|A^\circ|\|_{L^2(S^2)}
\|\e^{-\lambda}\delta\mathcal{W}\|_{L^2(S^2)}\\
&\leq C_p\sqrt{\mathcal{W}_0(\Phi)} \,
\|\e^{\lambda}\delta\mathcal{W}\|_{L^2(S^2)}\\
&\leq C_p\|\e^{\lambda}\delta\mathcal{W}\|_{L^2(S^2)},
\end{align}
where we estimated 
$\|\e^{-\lambda}\delta\mathcal{W}\|_{L^2(S^2)}
\leq \|\e^{\lambda}\delta\mathcal{W}\|_{L^2(S^2)}$
thanks to the simple inequality
\begin{align}
1=\e^{-\lambda}\e^{\lambda}\leq \sup_{S^2}(\e^{-\lambda})\e^{\lambda}\leq C(1+\delta)\e^{\lambda}
\end{align}
issuing from  property \eqref{eq:def-class-DLM}.
We now examine $\Omega$.
It will be more practical to look at the 
associated (time-dependent) conformal Killing vector field
i.e. generating 
conformal transformations:
\begin{align}
V=\Omega+\overline{\Omega},
\end{align}
Similarly as in the proof of Lemma \ref{lemma:key-IFT},
a basis for the vector space of
conformal Killing fields $T_e\Aut(S^2)$ is given,
in its immersed representative, is given by
\begin{align}
dI(Z_1)(y) &=
(0,-y^3,y^2), &
dI(Z_2)(y) &=
(y^3,0,-y^1),&
dI(Z_3)(y) &=
(-y^2,y^1,0),\\[0.5em]
dI(Z_4)(y) &=
e_1-y^1 y, &
dI(Z_5)(y) &=
e_2-y^2 y,&
dI(Z_6)(y) &=
e_3-y^3 y.&
\end{align}
One checks that this basis is orthogonal
with respect the $L^2$--scalar product
and each element has the same length.
Thus, we may write
\begin{align}
\label{eq:V-L2-basis}
V 
=\sum_{a=1}^6 V^a Z_a
=C\sum_{a=1}^6 (V,Z_a)_{L^2(S^2)} Z_a
=C\sum_{a=1}^6 (U,Z_a)_{L^2(S^2)} Z_a,
\end{align}
where the last inequality is a consequence of the fact that,
by construction, the normal solution
of the $\overline{\partial}$-operator is $L^2$-orthogonal
to the space of holomorphic vector fields.
So, to estimate $V$ it suffices to estimate
the (time-dependent) coefficients
\begin{align*}
(V,Z_a)_{L^2}
=\int_{S^2} \langle U, Z_a\rangle\,d\sigma
\quad\text{for }a=1,\ldots,6.
\end{align*}
Now, we may write the integrand as
\begin{align}
\langle U,Z_a\rangle
&=\langle dI(U),dI(Z_a)\rangle
=\langle d\Phi(U),dI(Z_a)\rangle
+\langle dI(U)-d\Phi(U),dI(Z_a)\rangle
\end{align}
where the second term can be estimated as
\begin{align}
\big|\langle dI(U)-d\Phi(U),dI(Z_a)\rangle\big|
\leq C|dI-d\Phi||U|
\end{align}
and thus, upon integration,
H\"older's inequality and
 property \eqref{eq:def-class-DLM}
of the set $\mathscr{W}^{\varepsilon,\delta}_{[0,T)}(S^2,\R^3)$,
we see that the estimate
\begin{align}
\label{eq:suffic-estimate}
\bigg|\int_{S^2}
\langle dI(U)-d\Phi(U),dI(Z_a)\rangle
d\sigma\bigg|
&\leq C\int_{S^2}|dI-d\Phi||U|\,d\sigma\\
&\leq C\|dI-d\Phi\|_{L^{p'}(S^2)}\|U\|_{L^p(S^2)}\\
&\leq C\|\Phi-I-c\|_{W^{2,2}(S^2)}\|U\|_{L^p(S^2)}\\
&\leq C_{p}\,\delta\, \|U\|_{L^p(S^2)},
\end{align}
holds, which will suit our purposes.
We are left to estimate the terms
\begin{align}
\int_{S^2}\langle d\Phi(U),dI(Z_a)\rangle\,d\sigma
\quad\text{for }a=1,\ldots,6.
\end{align}
Differentiating the well-balanced conditions 
\eqref{eq:well-balanced} gives,
\begin{align}
0&=\frac{d}{d t}
\int_{S^2} I\times\Phi\,d\sigma
=\int_{S^2}I\times(-\delta\mathcal{W} + d\Phi(U))\,d\sigma,\\
0&=\frac{d}{dt}\int_{S^2}
Id\sigma_g
=\int_{S^2}I
\big(\langle 2H,\delta\mathcal{W}\rangle
+\divop_g(U)\big)d\sigma_g,
\end{align}
and thus that
\begin{align}
\label{eq:diff-well-bal-1}
\int_{S^2}I\times d\Phi(U)\,d\sigma
&=\int_{S^2}I\times \delta\mathcal{W}\,d\sigma,
\end{align}
and,
upon integration by parts,
that
\begin{align}
\label{eq:diff-well-bal-2}
\int_{S^2}dI(U)\,d\sigma_g
&=\int_{S^2} I\langle 2H,\delta\mathcal{W}\rangle
\,d\sigma_g.
\end{align}
On the other hand,
a direct calculation shows that
\begin{align}
\label{eq:ptwise-scal-gen-1}
\begin{pmatrix}
\langle d\Phi(U), dI(Z_1)\rangle\\
\langle d\Phi(U), dI(Z_2)\rangle\\
\langle d\Phi(U), dI(Z_3)\rangle
\end{pmatrix}
=I\times d\Phi(U),
\end{align}
and similarly that
\begin{align}
\label{eq:ptwise-scal-gen-2}
\begin{pmatrix}
\langle d\Phi(U), dI(Z_4)\rangle\\
\langle d\Phi(U), dI(Z_5)\rangle\\
\langle d\Phi(U), dI(Z_6)\rangle
\end{pmatrix}
=d\Phi(U)-\langle d\Phi(U),I\rangle I
=(d\Phi(U))^\top,
\end{align}
where $(\cdot)^\top$ denotes the orthogonal projection
onto the tangent space of $S^2$ in the standard immersion.
Integrating \eqref{eq:ptwise-scal-gen-1}
and using \eqref{eq:diff-well-bal-1}
we can estimate for $a=1,2,3$
\begin{align}
\label{eq:est-coeff-123}
\big|(d\Phi(U), dI(Z_a))_{L^2}\big|
&\leq C\Big|\int_{S^2}I\times\delta\mathcal{W}\,d\sigma\Big|
\leq C\|\e^{\lambda}\delta\mathcal{W}\|_{L^2(S^2)}.
\end{align}
Integrating \eqref{eq:ptwise-scal-gen-2}
we get instead
\begin{align}
\sum_{a=1}^3
C(d\Phi(U),dI(Z_a))_{L^2} e_a
=\int_{S^2}(d\Phi(U))^\top d\sigma,
\end{align}
and if we write the integrand as
\begin{align}
(d\Phi(U))^\top
=(dI(U))^\top + (d\Phi(U)-dI(U))^\top
\end{align}
and notice that, similarly as for 
\eqref{eq:suffic-estimate} we can estimate
\begin{align}
\Big|\int_{S^2}(d\Phi(U)-dI(U))^\top d\sigma\Big|
&\leq C\int_{S^2} |d\Phi(U) - d I(U)|\,d\sigma\\
&\leq C\|d\Phi-dI\|_{L^{p'}(S^2)}\|U\|_{L^p(S^2)}\\
&\leq C_p\, \delta\,\|U\|_{L^p(S^2)},
\end{align}
using \eqref{eq:diff-well-bal-1} 
(since $dI(U) = dI(U)^\top$), we 
have for $a=4,5,6$
\begin{align}
\label{eq:est-coeff-456}
|(d\Phi(U),dI(Z_a))_{L^2}|
&\leq C \Big|\int_{S^2} d I(U)d\sigma\Big|
+C_p\, \delta\,\|U\|_{L^p(S^2)}\\
&\leq C\Big|
\int_{S^2} I\langle 2H,\delta\mathcal{W}\rangle
\,d\sigma_g\Big|
+C_p\, \delta\,\|U\|_{L^p(S^2)}\\
&\leq C\|H\e^{\lambda}\|_{L^2(S^2)}
\|\e^{\lambda}\delta\mathcal{W}\|_{L^2(S^2)}
+C_p\, \delta\,\|U\|_{L^p(S^2)}\\
&\leq C_p\big(\|\e^{\lambda}\delta\mathcal{W}\|_{L^2(S^2)}
+  \delta\,\|U\|_{L^p(S^2)}\big),
\end{align}
where we also used that 
$\|H\e^{\lambda}\|_{L^2(S^2)}=C\mathcal{W}_1(\Phi)$
is bounded uniformly in $t$.
\medskip

Estimates
\eqref{eq:suffic-estimate},
\eqref{eq:est-coeff-123} and \eqref{eq:est-coeff-456}
inserted in \eqref{eq:V-L2-basis}
yield
\begin{align}
\label{eq:est-V}
\|V\|_{L^\infty(S^2)}\leq C_p\big(\|\e^{\lambda}\delta\mathcal{W}\|_{L^2(S^2)}
+  \delta\,\|U\|_{L^p(S^2)}\big).
\end{align}
and so,
in conjunction with \eqref{eq:est-U-ellipt},
we get
\begin{align}
\|U\|_{L^p(S^2)}
\leq C_p
\big(\|\e^{\lambda}\delta\mathcal{W}\|_{L^2(S^2)}
+ \delta\,\|U\|_{L^p(S^2)}\big).
\end{align}
Taking the $L^2$-norm in time of such inequality
gives
\begin{align}
\|U\|_{L^2((0,T),L^p(S^2))}
\leq C_p
\big(\|\e^{\lambda}\delta\mathcal{W}\|_{L^2((0,T),L^2(S^2))}
+ \delta\,\|U\|_{L^2((0,T),L^p(S^2))}\big),
\end{align}
and so, if $\delta$ is chosen sufficiently small,
we reach
\eqref{eq:est-U}.
\end{bfproof}

\begin{remark}
\label{rmk:equiv-eps-delta}
The $\delta$
of Proposition \ref{prop:est-U}
may be smaller than that given by 
Proposition \ref{prop:well-balanced-fam}.
An inspection of the proof however 
shows however that we may equivalently have taken the same
$\delta$ at the price of choosing $\varepsilon>0$
sufficiently small,
because we may apply instead
Propositions \ref{prop:nonsmooth-area-DLM}
and \ref{prop:well-balanced-fam} for a.e. $t$ so that
$\Phi(t,\cdot)$ is $W^{4,2}$,
and argue in the end similarly as above.
This variant would have been
equally fine for our purposes.
\end{remark}

\begin{corollary}
\label{cor:bootstrap}
There exists $\varepsilon_0,\delta>0$ so that
any Willmore flow in
$\mathscr{W}^{\varepsilon_0,\delta}_{[0,T)}(S^2,\R^3)$
is in $C^\infty((0,T]\times S^2)$,
and, if it is a solution
to the Cauchy problem \eqref{eq:Cauchy-conformal-Willmore}
for smooth initial datum $\Phi_0$,
then it is in $C^\infty([0,T]\times S^2)$.
\end{corollary}

\begin{bfproof}
By definition and by Proposition \ref{prop:est-conform-flow},
$d\Phi$ is in
$L^\infty((0,T),W^{1,2}(S^2))\cap L^2((0,T),W^{3,2}(S^2))$,
we have sufficient regularity
to expand the Willmore operator in the flow equation
by means of formulas \eqref{eq:Willmore-divergence-form}
and \eqref{eq:bilaplace-decomp}:
\begin{align*}
\frac{\partial}{\partial t}\Phi
+\frac{1}{2}\e^{-4\lambda}\Delta^2\Phi
&=\frac{1}{2}\e^{-4\lambda}\big(
2\langle d\lambda,\nabla H\rangle_g 
+(2|d\lambda|^2_g + \Delta_g\lambda)H\big)\\
&\phantom{{}--{}}-\nabla^{*_g}\big(
\langle A, H\rangle^{\sharp_g}
+\langle A^\circ, H\rangle^{\sharp_g}
\big) + U;
\end{align*}
since $\e^{\lambda}$ is by assumption uniformly bounded,
the equation is uniformly parabolic,
and by Proposition \ref{prop:est-U},
$U\in L^2(0,T),L^p(S^2))$
for any $p<2$.

This is enough to start a boostrapping procedure
using first $L^p$-$L^q$ and then Schauder parabolic estimates
in a fashion similar to the elliptic case discussed in
Proposition \ref{prop:bootstrap-cond-Willmore}.
To bootstrap the regularity of the tangential component
$U$, one uses higher-regularity variants of Proposition
\ref{prop:est-U}, whose proofs are similar to the basic case.
\end{bfproof}

\begin{bfproof}[Proof of Theorem \ref{thm:main}]
\textit{Case of smooth initial datum.}
By Corollary \ref{cor:bootstrap}, 
it suffices to prove that there exists a unique 
smooth solution in $\mathscr{W}^{\varepsilon_0,\delta}_{[0,T)}(S^2,\R^3)$ with the
required properties.

An application  DeTurck's trick
(see Appendix \ref{sec:DeTurck}) yields existence
and uniqueness of a
smooth solution to the Cauchy problem for the
normal Willmore flow:
\begin{align}
\left\{
\begin{aligned}
\frac{\partial}{\partial t}\Phi^0 &= -\delta\mathcal{W}
&&\text{in }(0,T)\times S^2,\\
\Phi^0(0,\cdot)&=\Phi_0
&&\text{on } S^2,
\end{aligned}
\right.
\end{align}
and if $\varepsilon_0>0$ is small enough,
by Theorem \ref{thm:area-est}
$\Phi^0$ exists for all and smoothly converges
to a round sphere.%
\footnote{Note carefully: the convergence is
to \emph{some} smooth
parametrization, of a round sphere of 
\emph{some} center and radius,
not to $I$  modulo dilation and translation.
For instance, if $\Phi_0$ is \emph{any} smooth parametrization of 
a round sphere,
then $\Phi^0$ trivially converges to $\Phi_0$.}
We conformalize such flow composing it with the family
$(\phi(t,\cdot))_{t\in [0,+\infty)}$
of canonical quasi-conformal mappings associated
to the family of metrics $g^0(t,\cdot)=\Phi^0(t,\cdot)^\star g_{\R^3}$,
see  \cite{MR115006}.
The fact that it is $\phi(0,\cdot) = e$ that such family is smooth
both in the space and in time follows from the theory of quasi-conformal
mappings.
Then $\Phi^1(t,\cdot)=\Phi^0(t,\phi(t,\cdot))$ is a conformal Willmore flow
defined for all times
and converging to a conformal parametrization of a round sphere.
\medskip 

Let
\begin{align}
a(t) = \sqrt{\frac{4\pi}{\mathcal{A}(\mathcal{S}_t)}}
\end{align}
the normalizing function of time so that
$a(t)\mathcal{S}_t$ has always area $4\pi$
and let $\varepsilon_0>0$ be sufficiently small as in Theorem
\ref{thm:DeLellis-Muller} and also so that
$C\varepsilon_0\leq\delta$
where $C$ is as in \eqref{eq:DeLellis-Muller-est}
and $\delta$ is sufficiently small
as in Propositions  \ref{prop:well-balanced-fam} 
and \ref{prop:est-U}.%
\footnote{See in this regard Remark \ref{rmk:equiv-eps-delta}.}
In this way, there exists a family of conformal diffeomorphisms
$(\psi(t,\cdot))_{t\in[0,+\infty)}\subset \Aut(S^2)$
so that $\Phi(t,\cdot) = \Phi^1(t,\psi(t,\cdot))$ 
is conformal, well-balanced and
\begin{align}
\label{eq:prf-main-2}
\|a(t)\Phi(t,\cdot) - I-c(t)\|_{W^{2,2}(S^2)}
+\|a(t)\e^{\lambda(t,\cdot)}-1\|_{L^\infty(S^2)}
\leq C\sqrt{\mathcal{W}_0(\mathcal{S}_t)},
\end{align}
where $c(t)=\dashint_{S^2}\Phi(t,\cdot)\,d\sigma$
(see Remark \ref{rmk:minimality-avg})
and moreover there exist a neighborhood of the identity 
$\mathcal{O}\subset\Aut(S^2)$ where such choice is unique.

As for $a(t)$,
thanks the area control of Theorem \ref{thm:area-est} 
(recall that
$\mathcal{A}(\mathcal{S}) =  \mathcal{A}(\mathcal{S}_0)=4\pi$)
it is
\begin{align}
|\mathcal{A}(\mathcal{S}_t) - 4\pi|
\leq C\mathcal{W}_0(\mathcal{S}_t)
=o(1)
\quad\text{as } t\to+\infty,
\end{align}
and hence
\begin{align}
|a(t) - 1|
\leq C\mathcal{W}_0(\mathcal{S}_t)
=o(1)
\quad\text{as } t\to+\infty,
\end{align}
which means that we may remove $a(t)$ from the estimate
\eqref{eq:prf-main-2}.

As for $c(t)$, we may write
\begin{align}
c(t)
= \dashint_{S^2} \Phi(t,\cdot)\,d\sigma
&= \dashint_{S^2}
\Phi(t,\cdot)\e^{2\lambda(t,\cdot)}\,d\sigma
+\dashint_{S^2} \Phi(t,\cdot)(1-\e^{2\lambda(t,\cdot)})d\sigma\\
&=\mathcal{C}(\mathcal{S}_t) +
\dashint_{S^2} \Phi(t,\cdot)(1-\e^{2\lambda(t,\cdot)})d\sigma,
\end{align}
where $\mathcal{C}(\mathcal{S}_t)$ denotes the barycenter
as in Theorem \ref{thm:area-est}.
Thus, from:  the barycenter control of
Theorem \ref{thm:area-est},
the control on the conformal factor
issuing from \eqref{eq:prf-main-2},
the fact that  $\psi(t,\cdot)\in \mathcal{O}$,
smooth convergence
and $\mathcal{C}(\mathcal{S})	=0$,
we can estimate
\begin{align}
|c(t)|\leq C\mathcal{W}_0(\mathcal{S}_t),
\end{align}
and hence also $c(t)$ can be removed from estimate
\eqref{eq:prf-main-2} as well and deduce that
\begin{align}
\label{eq:prf-main-3}
\|\Phi(t,\cdot) - I\|_{W^{2,2}(S^2)}
+\|\e^{\lambda(t,\cdot)}-1\|_{L^\infty(S^2)}
\leq C\sqrt{\mathcal{W}_0(\mathcal{S}_t)}.
\end{align}
Since $\mathcal{W}_0(\mathcal{S}_t)=o(1)$
as $t \to+\infty$,
we obtain
that $\Phi(t,\cdot)$ converges to $I$ in
$W^{2,2}$
and that its conformal factor converges uniformly to $1$.
Convergence of higher order derivatives to $I$
then follow from this and the smooth convergence of $\Phi^1$.
\medskip

\textit{Case of weak initial datum.}
Let $\Phi_{0,j}\in \mathscr{D}^{\varepsilon_0}(S^2,\R^3)$
be a sequence appoximating $\Phi_0$
in the weak $W^{2,2}$-topology,
i.e.
\begin{align}
\label{eq:weak-suit-1}
\Phi_{0,j}&\rightharpoonup \Phi_0\quad\text{in }W^{2,2}(S^2).
\end{align}
If $\varepsilon_0$ is taken sufficiently small,
it follows from the analysis in
\cite{MR2430975,MR3276154,MR3524220}
(see for instance the proof of Theorem 3.36
in \cite{MR3524220}) that also
\begin{align}
\label{eq:weak-suit-2}
\e^{\lambda_j}&\overset{*}{\rightharpoonup}
\e^{\lambda}\quad\text{in }L^\infty(S^2).
\end{align}
For each $j$, we let $\Phi_j\in
\mathscr{W}^{\varepsilon_0,\delta}_{[0,+\infty)}(S^2,\R^3)$
be the well-balanced conformal Willmore flow 
given by Theorem \ref{thm:main}
with initial datum $\Phi_{0,j}$, which we know to be smooth.
By estimate
\eqref{eq:prf-main-3}
we have that for every $t>0$ there holds
\begin{align}
\|\Phi_j(t,\cdot) - I\|_{W^{2,2}(S^2)}
+\|\e^{\lambda_j(t,\cdot)}-1\|_{L^\infty(S^2)}
\leq C\sqrt{\mathcal{W}_0(\mathcal{S}_t)}
\leq \delta,
\end{align}
and by Proposition \ref{prop:est-conform-flow}
it follows that
for every fixed choice of $T>0$ there holds
\begin{align}
\label{eq:add-1}
\|\Phi_j\|_{L^2((0,T), W^{4,2}(S^2))}
\leq C,
\end{align}
and finally by Proposition \ref{prop:est-U}
we also have, for a fixed choice of $1<p<2$,
\begin{align}
\|\partial_t\Phi_j\|_{L^2((0,T),L^p(S^2))}
=\|-\delta\mathcal{W}_j + U_j\|_{L^2((0,1),L^p(S^2))}
\leq C,
\end{align}
where $C$ does not depend on $j$.
\medskip

Weak sequential compactness properties of Sobolev spaces then
imply that,
up to the extraction of a subsequence,
there exists measurable functions
${\Phi:(0,+\infty)\times S^2\to\R^3}$ and
${\lambda:(0,+\infty)\times S^2\to\R}$ so that,
for every fixed $T>0$, as $j\to\infty$, 
\begin{align}
\label{eq:add-2}
\Phi_j&\to\Phi\quad\text{a.e. in }(0,+\infty)\times S^2,\\
\label{eq:add-3}
\Phi_j&\rightharpoonup\Phi
\quad\text{in } W^{1,p}((0,T)\times S^2),\\
\Phi_j&\rightharpoonup\Phi
\quad\text{in } L^2((0,T), W^{4,2}(S^2)),\\
\Phi_j&\overset{*}{\rightharpoonup}\Phi
\quad\text{in } L^\infty((0,T), W^{2,2}(S^2)),\\
\label{eq:add-4}
\e^{\lambda_j}&\overset{*}{\rightharpoonup}\e^{\lambda}
\quad\text{in } L^\infty((0,T)\times S^2).
\end{align}
Thus, for every $T>0$,
 $\Phi$ satisfies the conditions
(i), (ii), (iv) of Definition \ref{def:well-balanced}
(and $\e^\lambda$ is its conformal factor),
its Willmore operator is in $L^2((0,T)\times S^2)$
and finally since also 
\begin{align}
\mathcal{W}_0(\Phi(t,\cdot))
\leq\liminf_{j\to\infty}\mathcal{W}_0(\Phi_j)
\leq\varepsilon_0,
\end{align}
which means that (ii) is satisfied as well
and so $\Phi\in \mathscr{W}^{\varepsilon_0,\delta}_{[0,+\infty)}(S^2,\R^3)$.
Moreover, since $p>1$
by \eqref{eq:add-3} $\Phi$ has a trace
at initial time, which, by uniqueness
and continuity of the trace operator, must coincide with 
$\Phi_0$.

Finally, the convergence properties 
\eqref{eq:add-2}-\eqref{eq:add-4}
are enough to pass
to the limit as $j\to\infty$
in Definition \ref{def:weak-Willmore-flow}
and thus deduce that $\Phi$ is also a weak Willmore flow.
By Corollary \ref{cor:bootstrap}
$\Phi$ is smooth on $(0,+\infty)\times S^2$ as well.
\end{bfproof}

\appendix

\section{DeTurck's Trick for 
the Willmore Flow}
\label{sec:DeTurck}
We outline here one way to obtain,
in the smooth category,
short-time existence for the Cauchy's problem
\eqref{eq:Cauchy-normal-Willmore}
adapting an idea originally devised by 
\textsc{DeTurck} \cite{MR697987} in the context
of the Ricci flow.
There are other possibilities,
such as  the graph Ansatz \cite{MR1731639}
(in codimension 1)
or through the Nash-Moser Implicit Function Theorem \cite{MR656198,MR664497}.

The key idea is as follows.
From the divergence form of the Willmore operator
\eqref{eq:Willmore-divergence-form}, we write
\begin{align}
\delta\mathcal{W}
&=\Delta_g H
+\nabla^{*_g}\big( 
\langle A^\circ,H\rangle^{\sharp_g}
+\langle A,H\rangle^{\sharp_g}\big)
=\frac{1}{2}\Delta_g^2\Phi
+F\left(
d\Phi,
\nabla^2\Phi,
\nabla^3\Phi
\right),
\end{align}
where $F$ is a smooth function
and the covariant derivatives 
are with respect to a fixed smooth reference
metric on $\Sigma$.
Now in local coordinates
we may write
\begin{align}
\Delta_g \Phi
&=g^{\mu\nu}\partial^2_{\mu\nu}\Phi
-g^{\mu\nu}\Gamma_{\mu\nu}^\sigma\partial_\sigma\Phi,
\end{align}
and so
\begin{align}
\Delta_g^2 \Phi
&=\Delta_g(\Delta_g\Phi)\\
&=g^{\mu\nu}\partial^2_{\mu\nu}(\Delta_g\Phi)
-g^{\mu\nu}\Gamma_{\mu\nu}^\sigma\partial_\sigma(\Delta_g\Phi)\\
&=g^{\mu\nu}\partial^2_{\mu\nu}
\left(g^{\alpha\beta}\partial^2_{\alpha\beta}\Phi
-g^{\alpha\beta}\Gamma_{\alpha\beta}^\gamma\partial_\gamma\Phi
\right)
+f\left(\partial_x\Phi,
\partial^2_{xx}\Phi,\partial^3_{xxx}\Phi\right),
\end{align}
for some smooth function $f$.
Now, because also the metric depends on $\Phi$:
$g_{\mu\nu} = \langle \partial_\mu\Phi,\partial_\nu\Phi\rangle$,
the term
$
g^{\alpha\beta}
\left(\partial^2_{\mu\nu}
\Gamma_{\alpha\beta}^\gamma
\right)
\partial_\gamma\Phi
$
also contains derivative of order 4 in $\Phi$.
This causes the Willmore operator to be degenerate elliptic,
and the corresponding flow to be degenerate parabolic.

However, writing:
\begin{align}
\Delta^2_{g}\Phi
&=g^{\mu\nu}g^{\alpha\beta}
\partial^4_{\mu\nu\alpha\beta}\Phi
-g^{\mu\nu}
\partial^2_{\mu\nu}\left(
g^{\alpha\beta}
\Gamma_{\alpha\beta}^\gamma\right)\partial_\gamma\Phi
+f\left(\partial_x\Phi,
\partial^2_{xx}\Phi,\partial^3_{xxx}\Phi\right),
\end{align}
one guesses that it may be possible to add a tangent vector field
to the Willmore operator which for which the corresponding flow is a 
uniformly parabolic one. This motivates the following.

\begin{definition}
\label{def:DeTurck-vect-field}
Let $\Sigma$ be a closed, orientable surface
and let $\Phi_0:\Sigma\to\R^3$ be a smooth immersion.
If $\Phi:\Sigma\to\R^3$ is another smooth immersion,
\emph{DeTurck's vector field for the Willmore flow} 
(for $\Phi$ relative to $\Phi_0$)
is the vector field tangent to $\Phi$ given by
\begin{align}
V = V(\Phi_0,\Phi)
=-\frac{1}{2}\Delta_g W
=-\frac{1}{2} (\Delta_g W)^\gamma \partial_\gamma\Phi,
\end{align}
where $W = W^\gamma\partial_\gamma$ is the vector field
on $\Sigma$ given by
\begin{align}
W^\gamma=
g^{\alpha\beta}
\left(
\Gamma_{\alpha\beta}^\gamma - \breve{\Gamma}_{\alpha\beta}^\gamma
\right),
\end{align}
where $\Gamma_{\alpha\beta}^\gamma$
and
$\breve{\Gamma}_{\alpha\beta}^\gamma$
denote respectively the Christoffel symbols of $\Phi$
and $\Phi_0$.
\end{definition}
Note that this definition makes sense, since it is a well-known fact in differential geometry that,
although the Christoffel symbols are themselves not tensor,
the expression $\Gamma_{\alpha\beta}^\gamma - \breve{\Gamma}_{\alpha\beta}^\gamma$
(and consequently its trace $W$) is.

\begin{proposition}
[Short-Time Existence for the Smooth DeTurck-Willmore Flow]
\label{prop:existence-Willmore-DeTurck-smooth}
Let $\Sigma$ be a closed, orientable surface
and let $\Phi_0:\Sigma\to\R^3$ be a smooth immersion.
There exists some
$T=T(\Phi_0)>0$ so that the Cauchy problem
\begin{align}
\label{eq:Willmore-DeTurck-Cauchy-pb}
\left\{
\begin{aligned}
\partial_t\Phi
&=-\delta\mathcal{W} + V
 &&\text{in }(0,T)\times \Sigma,\\
\Phi(0,\cdot)&=\Phi_0,
&&\text{on } \Sigma,
\end{aligned}
\right.
\end{align}
has a unique solution in the class $C^{\infty}([0,T]\times \Sigma,\R^3)$,
where  $V=V(\Phi_0,\Phi)$
is DeTurck's vector field for the Willmore flow.
\end{proposition}
\begin{bfproof}
It is sufficent to prove that \eqref{eq:Willmore-DeTurck-Cauchy-pb}
defines a uniformly parabolic system for $\Phi$ over $\Sigma$.
The existence, uniqueness and smoothness of a solution
follows then from the general theory for such systems
in H\"older spaces
\cite{MR0211083} (transl. English \cite{MR0211084}), \cite{MR0241822}.
We have
\begin{align}
\Delta_g  W = \tr_{g}\left(\nabla^{(2)} W\right)
&=g^{\mu\nu}\left(
\nabla^g_{\partial_\mu}\nabla^g_{\partial_\nu} W
- \nabla_{\nabla_{\partial_c}^g}^g W
\right)
=g^{\mu\nu}\left(
\nabla^g_{\partial_\mu}\nabla^g_{\partial_\nu} W
- \Gamma_{\mu\nu}^\sigma\nabla_{\partial_\sigma}^g W
\right),
\end{align}
so computing directly we see that
\begin{align}
\Delta_g W
&=g^{\mu\nu}\Big(
\partial^2_{\mu\nu}W^\xi
+\partial_\nu W^\sigma\Gamma_{\mu\sigma}^\xi
+\partial_\mu W^\sigma \Gamma_{\nu\sigma}^\xi\\
&\phantom{{}----{}}
-\Gamma_{\mu\nu}^\tau \partial_{\tau} W^\xi
+W^\sigma \partial_{\mu}\Gamma_{\nu\sigma}^\xi
+W^\sigma\Gamma_{\nu\sigma}^\tau \Gamma_{\mu\tau}^\xi
-\Gamma_{\mu\nu}^\sigma \Gamma_{\tau \sigma}^\xi W^\sigma
\Big)\partial_{\xi}\\
&=g^{\mu\nu}\left(
\partial^2_{\mu\nu}W^\xi
+f(W,\partial_x W,\partial_x\Phi,\partial^2_{xx}\Phi)
\right)\partial_\xi\\
&=g^{\mu\nu}\left(
\partial^2_{\mu\nu}W^\xi
+
f\left(\partial_x\Phi,
\partial^2_{xx}\Phi,\partial^3_{xxx}\Phi\right)
\right)\partial_\xi,
\end{align}
thus we have,
in \emph{every} choice of local coordinates,
that
\begin{align}
-\delta\mathcal{W} + V
&=
-\frac{1}{2}\Delta^2_g\Phi
+V
+f
\left(
\partial_x\Phi,
\partial^2_{xx}\Phi,\partial^3_{xxx}\Phi
\right)\\
&=
-\frac{1}{2}
\Big(
g^{\mu\nu}g^{\alpha\beta}
\partial^4_{\mu\nu\alpha\beta}\Phi
-g^{\mu\nu}
\partial^2_{\mu\nu}\left(
g^{\alpha\beta}
\Gamma_{\alpha\beta}^\gamma\right)\partial_\gamma\Phi\\
&\phantom{{}----{}}+g^{\mu\nu}\partial^2_{\mu\nu}
\left(
g^{\alpha\beta}
(\Gamma_{\alpha\beta}^\gamma - \breve{\Gamma}_{\alpha\beta}^\gamma)
\right)
\partial_\gamma\vec\Phi
\Big)
+f\left(\partial_x\Phi,
\partial^2_{xx}\Phi,\partial^3_{xxx}\Phi\right)\\
&=-\frac{1}{2}\left(
g^{\mu\nu}g^{\alpha\beta}
\partial^4_{\mu\nu\alpha\beta}\Phi
-g^{\mu\nu}\partial^2_{\mu\nu}
\left(
g^{\alpha\beta}
(\breve{\Gamma}_{\alpha\beta}^\gamma)
\right)
\partial_\gamma\vec\Phi
\right)
+f\left(\partial_x\Phi,
\partial^2_{xx}\Phi,\partial^3_{xxx}\Phi\right)\\
&=-\frac{1}{2}\left(
g^{\mu\nu}g^{\alpha\beta}
\partial^4_{\mu\nu\alpha\beta}\Phi
\right)
+f\left(\partial_x\Phi,
\partial^2_{xx}\Phi,\partial^3_{xxx}\Phi\right),
\end{align}
so that  \eqref{eq:Willmore-DeTurck-Cauchy-pb}
defines, for $T$ sufficiently small
and in the smooth category,
 a uniformly parabolic system of fourth order in $\Phi$.
\end{bfproof}
One now obtains the analogue short-time
existence and uniqueness statement for 
the Cauchy problem \eqref{eq:Cauchy-normal-Willmore}
combining Proposition \ref{prop:existence-Willmore-DeTurck-smooth}
and the the fact that there is a bijective correspondence
between tangential components an reparametrizations
(see the analogous discussion in \cite{MR2815949}).

\section{The $\overline{\partial}$-operator on Vector Fields}
\label{sec:delbar}
Let $\Sigma$ be a Riemann surface.
the complexified tangent bundle $T^{\C}\Sigma = T\Sigma + iT\Sigma$ splits in two sub-bundles:
\begin{align}
T^{\C}\Sigma = T\Sigma^{(1,0)} \oplus T\Sigma^{(0,1)},
\end{align}
whose sections are respectively $(1,0)$ and $(0,1)$-vector fields:
\begin{align}
\mathfrak{X}^{(1,0)}(\Sigma) &= \Gamma(T\Sigma^{(1,0)})
\quad\ni\quad V = (V^1+iV^2)\partial_z,\\
\mathfrak{X}^{(0,1)}(\Sigma) &= \Gamma(T\Sigma^{(0,1)})
\quad\ni\quad W = (W^1+iW^2)\partial_{\bar z}.
\end{align}
One can identify $(1,0)$- and real vector fields
by means of conjugation and $(1,0)$-projection:
\begin{align}
\mathfrak{X}^{(1,0)}(\Sigma)
\simeq \mathfrak{X}(\Sigma)
\quad : \quad
V\rightarrow V+ \overline{V}
\quad\text{and}\quad
 W^{(1,0)}\leftarrow W.
\end{align}
With this identification, we can consider holomorphic
vector fields as real vector fields,
whose flows consist precisely of 
families of $\Aut(\Sigma)$, the space of conformal self-maps
(M\"obius transformations) of $\Sigma$.
Calling $B=\Gamma(T\Sigma^{(1,0)}\otimes T^*\Sigma^{(0,1)})$,
the $\overline{\partial}$-operator over $(1,0)$-vector fields is
\begin{align}
\label{eq:delbar}
\overline{\partial}:\mathfrak{X}^{(1,0)}(\Sigma)\to 
B,
\quad
\overline{\partial}V = \frac{1}{2}\partial_{\bar{z}}
(V^1+iV^2)\partial_{z}\otimes d\bar{z},
\end{align}
whose kernel is the space of holomorphic vector fields
$
\ker(\overline{\partial}) = \mathfrak{X}^\omega(\Sigma).
$
Fix now  a conformal
metric over $\Sigma$, $g = \e^{2\lambda}|dz|^2$.
The formal $L^2$-adjoint of $\overline{\partial}$, 
defined through the formula
$(\overline{\partial}V, F)_{L^2} = (V, \overline{\partial}{}^*F)_{L^2}$ is,
if $F=f\partial_z\otimes d\bar{z}$,
\begin{align}
\label{eq:delbar-adj}
\overline{\partial}{}^*:B\to \mathfrak{X}^{(1,0)}(\Sigma),
\quad
\overline{\partial}{}^*F 
= - 2\e^{-4\lambda}\partial_z(\e^{2\lambda}f)\partial_z.
\end{align}
We deduce that $\ker(\overline{\partial}{}^*)$
consists of those tensors  $F=f\partial_z\otimes d\bar{z}$
so that $\e^{2\lambda}f$ is antiholomorphic.
If we lower the fist index of $F$:
\begin{align}
F{}^\flat = fg_{\bar{z} z}\,d\overline{z}\otimes d\overline{z}
=\frac{1}{2}e^{2\lambda}f\,d\overline{z}\otimes d\overline{z},
\end{align}
then $F^\flat$ is an antiholomorphic quadratic differential,
and so  its conjugate is a holomorphic quadratic differential.
With these identifications, we  have
\begin{align}
	\ker(\overline{\partial}{}^*)\simeq
Q^\omega(\Sigma).
\end{align}
In particular, we note that even though $\overline{\partial}{}^*$
does depend on the chosen metric,  $Q^\omega(\Sigma)$ does not.
\medskip

For given $F\in B$, we consider the equation
\begin{align}
\label{eq:delbar-pb}
\overline{\partial}V = F
\quad\text{on }\Sigma.
\end{align}
Then, \eqref{eq:delbar-pb} has a solution if and only if
\begin{align}
F\in \overline{\partial}(\mathfrak{X}^{(1,0)}(\Sigma))
=\ker(\overline{\partial}{}^*)^\perp
\end{align}
In such case, if $V_0$ is one such solution,
every other one is of the form 
$V= V_0+v$ for $v\in \mathfrak{X}^\omega(\Sigma)$.
\bigskip

The  \emph{normal solution} to \eqref{eq:delbar-pb} is 
the only one in $\ker(\overline{\partial})^\perp$,
and we denote it by $\overline{\partial}{}^{-1} F$.
Normal solutions satisfy the typical elliptic estimates,
such as for instance
\begin{align}
\|\overline{\partial}^{-1}F\|_{W^{1,2}(\Sigma)}\leq C \|F\|_{L^2(\Sigma)}
\quad\forall\, F\in\ker(\overline{\partial}{}^*)^\perp.
\end{align}
for a constant $C=C(\Sigma, g)>0$.
\medskip

As a consequence of the Riemann-Roch formula,
if $\gamma$ is the genus of $\Sigma$, we have
\begin{align}
\dim_{\C} Q^\omega(\Sigma)
=\begin{cases}
0 &\text{if }\gamma=0\\
1 &\text{if }\gamma=1,\\
3\gamma-3 &\text{if } \gamma\geq 2,
\end{cases}
\quad\text{and}\quad
\dim_{\C} \mathfrak{X}^\omega(\Sigma)
=
\begin{cases}
3 &\text{if }\gamma=0,\\
1 &\text{if } \gamma=1,\\
0 &\text{if } \gamma\geq 2.
\end{cases}
\end{align}
In particular,  \eqref{eq:delbar-pb} can be solved for any
$F$ when $\Sigma=S^2$.

{\small
\bibliographystyle{amsalpha}
\bibliography{WillmoreFlowSmallEnergyBibliography.bib}
}

\textsc{Francesco Palmurella}\\
Scuola Normale Superiore,
Piazza dei Cavalieri 7,
56126 Pisa, Italy\\
\texttt{francesco.palmurella@sns.it}
\bigskip

\textsc{Tristan Rivi\`ere}\\
ETH Z\"urich, D-MATH,
R\"amistrasse 101,
8092 Z\"urich, Switzerland\\
\texttt{tristan.riviere@math.ethz.ch}
\bigskip

\small
\textbf{Keywords}
Willmore surfaces, Willmore gradient flow,
Geometric flows,
Conformally invariant problems,
Conformal gauge, Weak solutions,
Critical non linear 4th order elliptic and parabolic equations.

\end{document}